\documentclass[12 pt]{article}
\usepackage{latexsym}
\usepackage{amsthm}
\usepackage{amssymb}
\usepackage{amsmath}
\usepackage{amsfonts}
\usepackage{mathrsfs}
\usepackage[all]{xy}
\usepackage{texdraw}
\usepackage{graphicx}
\usepackage{psfrag}
\usepackage{makeidx}

\newtheorem{Theory}{Theory}[section] %Counter for all types of Theorems
\newtheorem{theorem}[Theory]{Theorem}
\newtheorem{lemma}[Theory]{Lemma}
\newtheorem{technicalLemma}[Theory]{Technical Lemma}
\newtheorem{corollary}[Theory]{Corollary}
\newtheorem{proposition}[Theory]{Proposition}
\newtheorem{fact}{Fact}  %level of Exercise
\newtheorem{remark}[Theory]{Remark} %trivial but worth noticing
\newtheorem{question}{Question} %Open Question in my mind or in Theory
\newtheorem{conjecture}[question]{Conjecture}%Statement of suspicion

\newtheorem{Ntn}{Description} %Counter for all types of Notation/Definition
\newtheorem{Dn}[Ntn]{Definition}

\newcommand{\be}{\begin{enumerate}}
\newcommand{\ee}{\end{enumerate}}
\newcommand{\bq}{\begin{question}}
\newcommand{\eq}{\end{question}}
\newcommand{\bcj}{\begin{conjecture}}
\newcommand{\ecj}{\end{conjecture}}
\newcommand{\bc}{\begin{corollary}}
\newcommand{\ec}{\end{corollary}}
\newcommand{\bl}{\begin{lemma}}
\newcommand{\el}{\end{lemma}}
\newcommand{\btl}{\begin{technicalLemma}}
\newcommand{\etl}{\end{technicalLemma}}
\newcommand{\bt}{\begin{theorem}}
\newcommand{\et}{\end{theorem}}
\newcommand{\bp}{\begin{proposition}}
\newcommand{\ep}{\end{proposition}}
\newcommand{\bft}{\begin{fact}}
\newcommand{\eft}{\end{fact}}
\newcommand{\brk}{\begin{remark}}
\newcommand{\erk}{\end{remark}}
\newcommand{\bd}{\begin{Dn}}
\newcommand{\ed}{\end{Dn}}

\newcommand{\ploi}{PL_o(I) }

\newcommand{\Z}{ \mathbf Z }
 
\newcommand{\N}{\mathbf N }
 
\newcommand{\R}{\mathbf R }

\newcommand{\supp}[1]{\mathop{\textrm{Supp}}(#1)}

\newcommand{\brkp}[1]{\mathop{\mathscr{B}}(#1)}

\author{Collin Bleak}

\begin{document}
\centerline{{\Large{\bf A Minimal Non-solvable Group of
      Homeomorphisms}}}
\vspace{.1 in}

\centerline{Collin Bleak}
\centerline{Department of Mathematics, Binghamton University, Binghamton, NY 13902-6000 USA}
\vspace{.3 in}

\centerline{{\bf Abstract}}
\vspace{.1 in}

Let $PL_0(I)$ represent the group of orientation-preserving
piecewise-linear homeomorphisms of the unit interval which admit
finitely many breaks in slope, under the operation of composition.  We
find a non-solvable group $W$ and show that $W$ embeds in every
non-solvable subgroup of $PL_0(I)$.  We find mild conditions under
which other non-solvable subgroups ($B$, $(\wr\Z\wr)^{\infty}$,
$(\Z\wr)^{\infty}$, and $^{\infty}(\wr\Z)$) embed in
subgroups of $PL_0(I)$. We show that all solvable
subgroups of $PL_0(I)$ embed in all non-solvable subgroups of
$PL_0(I)$.  These results continue to apply if we replace
$PL_0(I)$ by any generalized R. Thompson group $F_n$.

\tableofcontents
\section{Introduction}
Let the symbol $\ploi$ represent the group of orientation-preserving
piecewise-linear homeomorphisms of the unit interval which admit
finitely many breaks in slope, under the operation of composition.  We
show that there is a non-solvable group $W$ so that $W$ embeds in
every non-solvable subgroup of $\ploi$.  We show that every solvable
subgroup of $\ploi$ embeds in every non-solvable subgroup of $\ploi$
(see \cite{bpgsc} and \cite{bpasc} for a geometric and two algebraic
classifications of the solvable subgroups of $\ploi$).  We show that
all virtually solvable subgroups of $\ploi$ are in fact solvable.
Finally, if $H\leq\ploi$, we find various mild conditions on the
action of $H$ on the unit interval which imply the existence in $H$ of
embedded copies of various of the non-solvable groups $B$,
$(\wr\Z\wr)^{\infty}$, $(\wr\Z)^{\infty}$ or $(\Z\wr)^{\infty}$ (these
last groups are defined in sections 1.2.3 and 1.2.5).

  We note that $\ploi$ has received attention from various
researchers lately, primarily because it is ``a source of groups with
interesting properties in which calculations are practical'' (quoting
Brin and Squier in \cite{picric}).  For example, $\ploi$ contains
copies of each of the generalized R. Thompson groups $F_n$, which have
themselves been a focus of current research.  All of our stated results
for subgroups of $\ploi$ hold if we replace the group $\ploi$ in each
of the statements with any particular group $F_n$.  (The groups $F_n$
were introduced by Brown in \cite{BrownFinite}, where they were
denoted $F_{n,\infty}$. These groups were later extensively studied by
Stein in \cite{SteinPLGroups}, by Brin and Guzm\'an in
\cite{BGAutomorphisms} and by Burillo, Cleary, and Stein in
\cite{BCS}.)

This paper is a logical continuation of the investigations began in
\cite{bpgsc} and continued in \cite{bpasc}, and it indirectly uses some
machinery developed in \cite{BrinEG}.  While we will only need a small
part of the theory developed in \cite{bpasc}, we will need almost all
of the definitions, theory, and techniques developed in \cite{bpgsc}.
Instead of simply restating the whole of \cite{bpgsc}, we will assume
the reader is familiar with that paper, although we will restate
relevant definitions here.

Our purely algebraic results are stated in section 1.1, while our
geometric results are stated in section 1.2, after we give the
required definitions.  We note in passing that the proof of our main
algebraic result depends in large part on our geometric results.

\subsection{Statements of Algebraic Results and Some History\label{algDefs}}
We now define some groups needed for a more precise statement of the
main results.  In order to do this, we must first recall the
definition of a standard restricted wreath product of groups.

Let $C$ and $T$ be groups.  Let $M = \bigoplus_{t\in T}C$ represent
the direct sum of copies of $C$ indexed by the elements of $T$ (as
opposed to the direct product).  We will
denote the group $M\rtimes T$ (where the action of $T$ on $M$ is by
right multiplication on the indices) by the symbol $C\wr T$.  The
group $C\wr T$ is the standard restricted wreath product of groups.
Following standard convention, we will refer to $C$ as the ``Bottom
group'' of $C\wr T$, we will refer to $M$ as the ``Base group'' of
$C\wr T$, and finally, we will refer to $T$ as the ``Top group'' of
$C\wr T$.  As we will not have a need to explicitly discuss other
types of wreath products in this paper, we will use the phrase
``Wreath product'' to mean the ``Standard restricted wreath product''
in the remainder.  Note that we can think of $C$, $M$, and $T$ as
subgroups of $C\wr T$ in fairly obvious manners (a little care is
required when choosing the realization of the group $C$ in $C\wr T$,
as many candidate copies of $C$ are available).

  If $C$ is a subgroup of $\ploi$, then there is a straightforward
geometric construction that realizes the group $C\wr\Z$ in $\ploi$.
We will give a concrete demonstration of this construction below in subsection
\ref{geoDefs}.  The construction is so basic in $\ploi$ that it
motivates the definition of the following family of groups, which will
play a central role in all that follows.

 Define
\[
\begin{array}{lcr}
W_0 = \left\{1\right\},& &\textrm{ and} \\
W_n = W_{n-1}\wr\Z,& &\textrm{ for $n>0$ with }n\in\N.
\end{array}
\] 

Note here that we use $\N = \left\{0,1,2,\ldots\right\}$.

We are now ready to define the minimal non-solvable subgroup in $\ploi$
and state our chief result.

First, define

\[
W = \bigoplus_{n\in\N}W_n.
\]

We have the following purely algebraic result.

\bt
\label{NSC} 

Let $H$ be a subgroup of $\ploi$.  $H$ is non-solvable if and only if
$W$ embeds in $H$.
\et

Let us place this result in the context of previous investigations.
We need to mention three non-solvable groups.  These groups are
$^{\infty}(\wr\Z)$, $(\Z\wr)^{\infty}$, and $B$. The first two of
these groups are direct limits of groups from the family
$\{W_n\}^{\infty}_{n = 1}$.  We will refer to the pair as the ``Limit
groups.''  The group $B$ contains embedded copies of the limit groups,
and is finitely generated. Brin defines all three groups in
\cite{BrinEG}, where he denotes $B$ as $G_1$ in section five of his
text.  (We will define these groups as well, in section \ref{geoDefs}
below.)

Sapir had asked the question of whether every non-solvable subgroup of
$F$ contained a copy of $(\Z\wr)^{\infty}$.  Brin in \cite{BrinEG}
answers this negatively by showing that both $(\Z\wr)^{\infty}$ and
$^{\infty}(\wr\Z)$ occur as embedded subgroups of $F$, and that
neither of these two groups contains the other as an embedded
subgroup.  Brin then asks (Question 1 of his text) whether one of these two
groups has to occur as an embedded subgroup in any non-solvable
subgroup of $\ploi$.  Since $W$ contains neither of these groups as
embedded subgroups, Theorem \ref{NSC} answers Brin's question in the
negative.  

Theorem \ref{NSC} has various other consequences.

It is a consequence of the main results of \cite{bpasc} that each
solvable subgroup $H$ of $\ploi$ admits an $n\in\N$ so that $H$ embeds
in $W_n$.  However, each $W_n$ embeds in $W$.  Therefore, we
have the following corollary.

\bc
\label{solveEmbeddings}

Every solvable subgroup of $\ploi$ embeds in every non-solvable
subgroup of $\ploi$.

\ec

Since $W$ is not virtually solvable, we see that the non-solvable
subgroups of $\ploi$ are not virtually solvable, in particular, we have
the following second corollary.

\bc
\label{VSC}

Suppose $H$ is a subgroup of $\ploi$.  If $H$ is virtually solvable,
then $H$ is solvable.

\ec

As mentioned above, the papers \cite{bpgsc} and \cite{bpasc} provide
both geometric and algebraic classifications of the solvable subgroups
of $\ploi$.  Therefore, the last corollary, together with the cited
research, is sufficient to classify the virtually solvable subgroups
of $\ploi$.

\subsection{Definitions and Statements of Geometric Results \label{geoDefs}}
We show stronger results than the above stated Theorem \ref{NSC},
under mild restrictions on the nature of the non-solvable subgroups of
$\ploi$ under investigation.  The extra hypotheses involved all have
a visual context, so that we will refer to all statements of results
in this section as ``Geometric results.''

Suppose $H$ is a subgroup of $\ploi$.  In subsection 1.2.3 we give a
geometric criterion which when satisfied implies that $H$ will admit a
copy of Brin's group $B$ as a subgroup. In subsection 1.2.5 we state
various geometric conditions which when individually satisfied imply
that $H$ will admit an embedded copy of a particular group in the list
$(\wr\Z\wr)^{\infty}$, $(\Z\wr)^{\infty}$, and $(\wr\Z)^{\infty}$.
The definitions of the geometric criteria mentioned above and of the
groups $(\wr\Z\wr)^{\infty}$, $(\Z\wr)^{\infty}$, $(\wr\Z)^{\infty}$,
and of $B$ are spread out through section 1.2.
 
For the next two subsections, we fix a model group $G\leq\ploi$ in
order to have a common reference while we define some terminology.

\subsubsection{General geometric definitions}

Define $\supp{G}$, the \emph{support of $G$}, to be the set
$\left\{x\in I\,|\, xg\neq x \textrm{ for some } g\in G\right\}$.  The
set $\supp{G}$ is an open subset of $(0,1)$, and can therefore be
written as a disjoint union of a countable (possibly finite)
collection of open intervals in $(0,1)$.  If $g\in G$ then we will
similarly refer to $\supp{\langle g\rangle}$ as the \emph{support of
$g$}.  We note in passing that if $g$, $h\in G$ and the support of $g$
and the support of $h$ are disjoint, then $h$ and $g$ commute.  We
call any interval component of $\supp{G}$ an \emph{orbital of $G$}.
If $g\in G$, and $A$ is an orbital of $\langle g\rangle$, then we say
$A$ is an \emph{orbital of $g$} or an \emph{element-orbital of $G$}.

We also note that if $g\in G$ and $A$ is an orbital of $g$ then either
$g$ moves all points in $A$ to the right or $g$ moves all points in
$A$ to the left.  If $g$ moves all points in $A$ to the right then we
will refer to the interval $[x,\,xg)$ as a \emph{fundamental domain of
$g$ in $A$} (note that in this paper, all group actions will be right
actions, and also that we will compose elements from the left to the
right).  If $g$ moves all points in $A$ to the left, then we will
similarly refer to the interval $(x,\,xg^{-1}]$ as a \emph{fundamental
domain of $g$ in $A$}.  We will occasionally not mention the orbital
$A$ if the context will allow us to do this without confusion.  In any
case, note that a fundamental domain of an element $g\in G$ in one of
its orbitals $A$ is a maximal subinterval of $A$ that is entirely
mapped off of itself by the action of $g$.

If $A$ is an element-orbital of $G$, there will be infinitely many
elements in $G$ with orbital $A$, so that we often explicitly
associate an element with an element-orbital.  We call an ordered pair
$(A,\,g)$ a \emph{signed orbital of $G$} if $g$ is an element of $G$
with orbital $A$.  In this case, we refer to $A$ as the \emph{orbital
of $(A,\,g)$} and $g$ as the \emph{signature of $(A,\,g)$}.  We will
often work with sets of signed orbitals.  If $X$ is a set of signed
orbitals, then we will sometimes form sets $O_X$ and $S_X$, where $O_X
= \left\{A\subset I\,|\,(A,\,g)\in X \textrm{ for some } g\in
G\right\}$ and $S_X = \left\{g\in G\,|\, (A,\,g)\in X \textrm{ for
some} A \subset I\right\}$.  We will refer to the set $O_X$ as
\emph{the orbitals of $X$} and to the set $S_X$ as \emph{the
signatures of $X$}.

\subsubsection{Transition Chains}

Suppose $\mathscr{C} = \left\{(A_p,g_p)\,|\,p\in \mathscr{I}\right\}$
is a set of signed orbitals indexed by a set $\mathscr{I}\subset I$
and $A_{\mathscr{C}} = \cup_{p\in\mathscr{I}}A_p$.  We call $\mathscr{C}$
a \emph{transition chain} if $\mathscr{C}$ satisfies the following
conditions:

\be
\item For all $x$, $y\in A_{\mathscr{C}}$, with $x<y$, the interval
$[x,y]\subset A_{\mathscr{C}}$.
\item For all $p\in \mathscr{I}$, if $g\in S_{\mathscr{C}}$ and
$p\in\supp{g}$, then $g = g_p$.  
\ee

In this case, we refer to the cardinality of $\mathscr{C}$ as the
\emph{length of $\mathscr{C}$}.  

The point of a transition chain is that the union $Z$ of the orbitals
of a transition chain represents an interval in $I$ on which the
signatures of the transition chain may act non-trivially.  In
particular, by carefully choosing which signatures to act with, and in
some specific order, we can move any particular point in $B$ as far to
the left or right in $Z$ as desired.

In the case of short transition chains, we will have the second
condition fairly easily from direct considerations.  For instance, here
is an alternative specific definition of a transition chain of length
two (we will not need any transition chains in this paper of length
greater than three), note that the second condition of the general
definition is satisfied here.  Suppose $\mathscr{C} =
\left\{(A,g),\,(B,h)\right\}$ is a set of signed orbitals of $G$,
where $A = (a_l,a_r)$ and $B = (b_l,b_r)$ and $a_l<b_l<a_r<b_r$.  In
this case we call $\mathscr{C}$ a \emph{transition chain of length two
for $G$}.

We now release the group $G$.  We will use the language developed
above freely with other subgroups of $\ploi$, expecting that this will
lead to no confusion. 

\subsubsection{The group $B$, and our chief geometric result}

With transition chains defined, we only need a definition of the group
$B$ in order to state our first geometric result.  The group $B$ was
introduced in a general form in \cite{BrinEG} under the notation
$G(5)$ in section 5 of that paper.  Our realization will be much more
concrete, but as in \cite{BrinEG}, it will still be realized in
R. Thompson's group $F$ (the subgroup of $\ploi$ consisting of
elements whose breaks in slope occur only in the dyadic rationals, and
which have all slope values powers of $2$).

Define $\alpha\in\ploi$ to be the element so that given any $x\in I$,
we have
\[
x\alpha = \left\{\begin{array}{ll}
\frac{1}{4}x & 0\leq x < \frac{1}{4},
\\
\\
x-\frac{3}{16} & \frac{1}{4} \leq x<\frac{7}{16},
\\
\\
4x - \frac{3}{2} & \frac{7}{16}\leq x<\frac{9}{16},
\\
\\
x+\frac{3}{16} & \frac{9}{16}\leq x<\frac{3}{4},
\\
\\
\frac{1}{4}x + \frac{3}{4}& \frac{3}{4}\leq x\leq 1,
\end{array}\right.
\]
and define $\beta_0\in\ploi$ to be the element so that given any $x\in
I$, we have
\[
x\beta_0 =\left\{\begin{array}{ll}
x&0\leq x<\frac{7}{16},
\\
\\
2x-\frac{7}{16}&\frac{7}{16} \leq x< \frac{15}{32},
\\
\\
x+\frac{1}{32} & \frac{15}{32}\leq x<\frac{1}{2},
\\
\\
\frac{1}{2}x + \frac{9}{32} & \frac{1}{2} \leq x < \frac{9}{16},
\\
\\
x & \frac{9}{16}\leq x\leq 1.
\end{array}\right.
\]
The graphs of these elements (superimposed) are given below.  

\begin{center}
\psfrag{a}[c]{$\alpha$}
\psfrag{b1}[c]{$\beta_0$}
\includegraphics[height=340pt,width=340 pt]{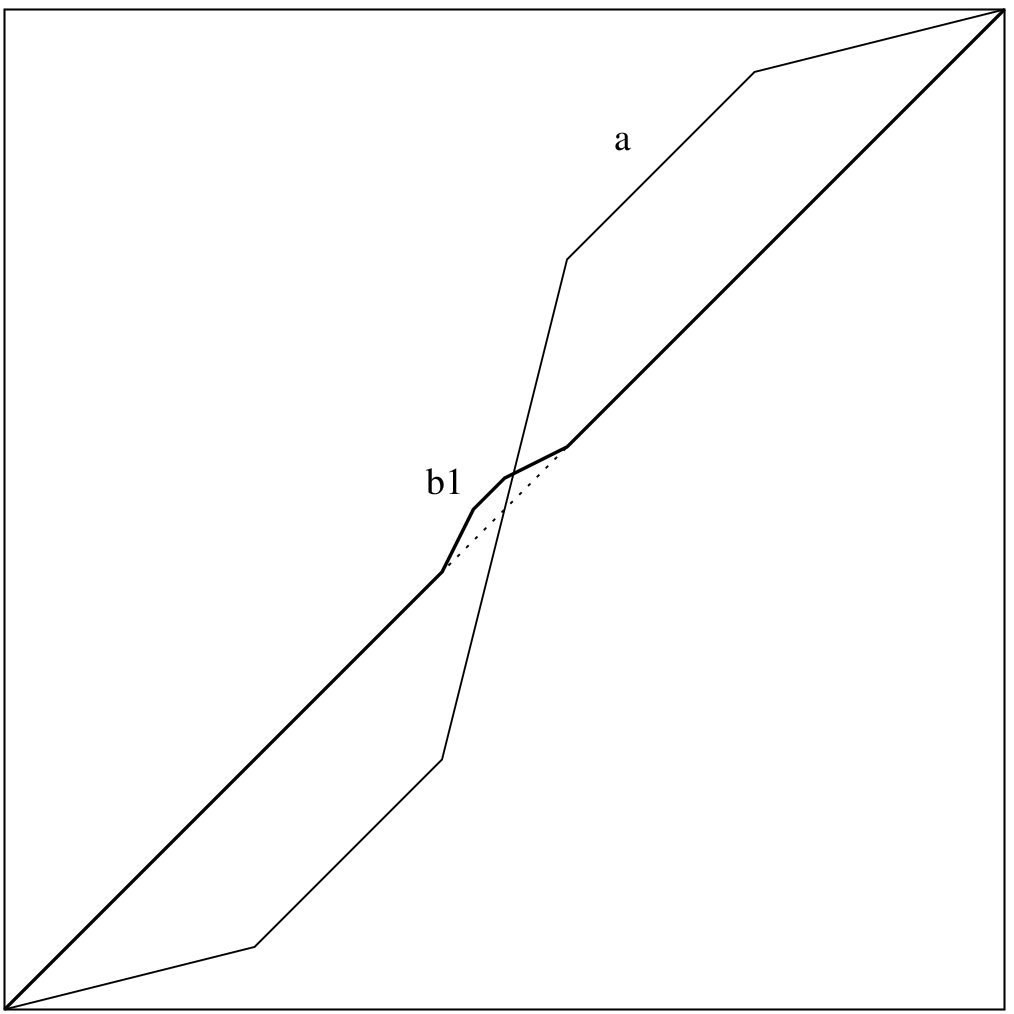}
\end{center}

We define
\[
B = \langle\alpha,\beta_0\rangle.
\]

We show, via a strengthening of an argument in
\cite{bpgsc}, the following result.

\bt
\label{tChainB}
If $H\leq \ploi$ admits a transition chain of length two, then $B$
embeds in $H$.
\et

It is extremely restrictive to only consider subgroups of $\ploi$
which do not admit transition chains of length two.  In particular,
$B$ embeds in ``most'' naturally occurring subgroups of $\ploi$.

\subsubsection{The structure of the group $B$}

We would like to discuss $B$ further, and extend our chief geometric
result, but we need to establish a few more conventions before we
proceed.  Let $x\in I$ and $g$, $h\in\ploi$, and recall that all group
actions will be written as right actions (so that $xg$ is the point
$x$ goes to when acted upon by $g$) and that compositions will occur
from the left to the right.  We will represent the conjugate of $g$ by
$h$ by the expression $g^h$, which will mean $h^{-1}gh$, and the
commutator of $g$ and $h$ by the expression $[g,h]$, which will mean
$g^{-1}h^{-1}gh$.

Let us now discuss the group $B$.  One perspective on the group $B$ is
that the element $\alpha$ acts as a ``Growing conjugator,''
conjugating $\beta_0$ to new elements
$\left\{\beta_i\right\}_{i\in\Z}$ (which elements, for positive index
$i$, have larger supporting sets).  We demonstrate this behavior, and
explore the structure of the groups generated by sub-collections of the
$\beta_i$.

Define $\beta_k =
\beta_0^{\alpha^k}$ for each integer $k$.  In particular, $\beta_{1}$
is given by the rule
\[
x\beta_1=x\alpha^{-1}\beta_0\alpha = \left\{\begin{array}{ll}
x &0\leq x<\frac{1}{4},
\\
\\
2x-\frac{1}{4} &\frac{1}{4}\leq x < \frac{3}{8},
\\
\\
x+\frac{1}{8} &\frac{3}{8}\leq x < \frac{1}{2},
\\
\\
\frac{1}{2}x+\frac{3}{8} & \frac{1}{2}\leq x<\frac{3}{4},
\\
\\
x&\frac{3}{4}\leq x\leq 1.
\end{array}\right.
\]
Observe that the support of $\beta_0$ is
$(\frac{7}{16},\frac{9}{16})$, and since $\frac{7}{16}\beta_1 =
\frac{9}{16}$, the support of $\beta_0$ is contained in a single
fundamental domain of $\beta_1$.  In particular, given any $i\in\Z$,
we see that the support of $\beta_{i-1}$ is contained in a single
fundamental domain of the support of $\beta_i$, since these two
elements are conjugates of $\beta_0$ and $\beta_1$.  Therefore, given
an $i\in\Z$, any two conjugates of $\beta_i$ by distinct powers of
$\beta_{i+1}$ will have disjoint support, so that these two conjugates
will commute.  In particular, the group $\langle
\beta_i,\,\beta_{i+1}\rangle$ is isomorphic with $\Z\wr\Z$.  More
generally, given an $n\in\N$, we see that any collection of $n$
distinct $\beta_i$ will generate a group isomorphic with $W_n$.

Let $\mathscr{W} = \left\{\beta_i\right\}_{i\in\Z}$ throughout the paper.

In passing, let us point out that if $G$ is a subgroup of $\ploi$,
then by considering $G$ to be a subgroup of $Homeo(\R)$ (extend the
elements of $G$ by using the identity function outside of the unit
interval), we can conjugate $G$ (by a conjugator in $Homeo(\R)$) to
$\hat{G}$, a piecewise linear copy of $G$ whose support is in
$(\frac{7}{16},\,\frac{9}{16})$.  By the discussion above, we see that
the group $\langle \hat{G}, \beta_0\rangle$ is therefore isomorphic to
$G\wr\Z$.  This is a concrete version of the construction mentioned
earlier in the introduction.

We can think of $B$ as an HNN extension of the group generated by the
full collection $\mathscr{W}$ of the $\beta_i$, where $\alpha$ plays
the role of the stable letter, with the rule that $\beta_i^{\alpha} =
\beta_{i+1}$ for each integer $i$.  This then gives us a criterion for
detecting when a two-generator (say $\omega_0$ and $\gamma$) subgroup
of $\ploi$ is isomorphic to $B$.  For each integer $i$, define
$\omega_i = \omega_0^{\gamma^{i}}$.  The group $\langle
\omega_0,\gamma\rangle$ will be isomorphic to $B$ if the set map
$\Upsilon:\left\{\omega_i\right\}_{i\in\Z}\to\left\{\beta_i\right\}_{i\in\Z}$,
defined by the rule $\omega_i\mapsto\beta_i$ for all integers $i$, is
well defined and extends to an isomorphism of the groups
$\langle\left\{\omega_i\right\}_{i\in\Z}\rangle\cong\langle\mathscr{W}\rangle$.
Note that this last isomorphism is easy to detect geometrically via
the tools developed in \cite{BrinEG}; in our case, if the closure of
the support of the non-trivial $\omega_{0}$ is fully contained in a
set $X_{0}$, and $X_{0}\omega_1\cap X_{0} = \emptyset$, then
$\langle\left\{\omega_i\right\}_{i\in\Z}\rangle\cong\langle\mathscr{W}$.

\subsubsection{Towers, and further results}
Our last set of geometric results is mostly relevant in the situation
where a subgroup of $\ploi$ fails to admit transition chains of length
two.  In this restrictive case, we use other sets of signed orbitals
to understand our group structures.

It is commonly known that $\ploi$ is totally ordered, and that the set
of open intervals in $I$ is a poset under inclusion.  Using the
induced lexical ordering, we see that the set of signed orbitals of
$\ploi$ is itself a poset.  We will use this to help form the following
definition.

Given a subgroup $G$ of $\ploi$, we say
that a set $\mathscr{T}$ of signed orbitals of $G$ is a \emph{tower
associated with $G$} if $\mathscr{T}$ satisfies two properties.

\be
\item $\mathscr{T}$ is a chain in the poset of signed orbitals of $\ploi$.
\item If $A\subset I$ and $(A,g)$ and $(A,h)$ are in $\mathscr{T}$,
then $g = h$.  
\ee

Note in passing that the second condition above assures us that
orbitals in a tower get larger as we move to larger elements in the
tower.

We define the \emph{height} of a tower to be its cardinality.  We
define the \emph{depth} of a subgroup $G$ of $\ploi$ to be the supremum of the
set of cardinalities of towers associated with $G$.

The following is a variation of the main result of
\cite{bpgsc}.  \bt
\label{GSC}

Suppose $H$ is a subgroup of $\ploi$.  $H$ is non-solvable if and only
if $H$ admits towers of arbitrary finite height.  
\et

This last result assures us that if $H\leq \ploi$ admits an infinite
tower, then $H$ will be non-solvable.  We now focus on the sorts
of subgroups we can find in such a group $H$, depending on the types
of transition chains and towers we can find associated with $H$.

If a tower $T$ admits an order preserving injection $t:\N\to T$, then
we will say that the tower is \emph{tall}.  If a tower $T$ admits an
order reversing injection $d:\N\to T$, then we will say that a tower
is \emph{deep}.  If a tower $T$ admits an order preserving injection
$b:\Z\to T$ then we will say $T$ is \emph{bi-infinite}.

We are about ready to state our last set of results; we first need
concrete realizations of the groups $^{\infty}(\wr\Z)$ and
$(\wr\Z)^{\infty}$, as well as another group $(\wr\Z\wr)^{\infty}$ (all
of which are defined in \cite{BrinEG}, in a less concrete fashion).

Define the following groups.
\[
\begin{array}{lcl}
(\wr\Z\wr)^{\infty}& = &\langle\mathscr{W}\rangle \\ 
&&\\ 
^{\infty}(\wr\Z)&= &\langle \left\{\beta_i\,|\,i<0, i\in\Z\right\}\rangle \\ 
&&\\
(\Z\wr)^{\infty}&= &\langle \left\{\beta_i\,|\,i\in\N\right\}\rangle \\
\end{array}
\]
Recalling that, for instance, $W_3\cong \langle \beta_{-3},
\beta_{-2}, \beta_{-1}\rangle$ while $W_2 \cong\langle\beta_{-2},
\beta_{-1}\rangle$, it is natural to think of the groups
$^{\infty}(\wr\Z)$ and $(\Z\wr)^{\infty}$ as limit groups built by
using different families of inclusion maps $W_i\to W_{i+1}$.

We are now ready to state our further results.

\bt
\label{tallWreath}
\label{deepWreath}
\label{biWreath}

Suppose $G$ is a subgroup of $\ploi$, then
\be

\item if $G$ admits a tall tower then $G$ contains a
subgroup isomorphic to $(\Z\wr)^{\infty}$,
\item if $G$ admits a deep tower then $G$ contains a subgroup
  isomorphic to $^{\infty}(\wr\Z)$, and
\item if $G$ admits a bi-infinite tower then $G$ contains a subgroup
isomorphic to $(\wr\Z\wr)^{\infty}$.
\ee 

\et

The above theorem may not seem surprising, given the realizations of
$^{\infty}(\wr\Z)$, $(\Z\wr)^{\infty}$, and $(\wr\Z\wr)^{\infty}$
above.  However, the collections of signatures from the towers in the
theorem might admit other orbitals, away from the specified towers, so
that the groups generated by the collections of the signatures of the
towers can exhibit much more complicated behavior then can be found in
$(\wr\Z\wr)^{\infty}$.  Removing this ``External complexity'' is the
main work of this paper.

There has been some work towards further strengthening Theorem
\ref{NSC} in the case of a finitely generated non-solvable subgroup of
$\ploi$.  The current status (see \cite{BDiss}) is that one of
$^{\infty}(\wr\Z)$ or $(\Z\wr)^{\infty}$ must embed as a subgroup of
$\ploi$.  None of that work will appear in this paper.

The author would like to thank Matt Brin and Binghamton University for
their support during the research leading up to this paper.  Some of
the results here are contained in the author's dissertation written at
Binghamton University.

\section{Some Essential Geometry}
In this section, we will review the necessary known geometric facts
about $\ploi$.  We will prove only one of the results in this section, as 
the remainder can be found in, or are straightforward consequences of,
the results in \cite{BSPLR, BrinU, bpgsc} and \cite{bpasc}.  (While
the result we prove is new, its proof is straightforward using the
ideas in \cite{bpgsc}, so we include it in this section.)  For the
more complex known results, we will indicate references more
precisely.

We begin with some often-used facts about the action of
elements of $\ploi$ on the unit interval.

\brk If $g\in\ploi$, $A$ is an orbital of $g$ then the
following are true.  \be
\item $g$ has finitely many orbitals.
\item Either $xg>x$ for all $x\in A$ or $xg<x$ for all $x\in A$.
\ee 

\erk 

Because of the second point above, given an orbital $A$ of an element
$g\in\ploi$, we will say that \emph{$g$ moves points right (left) on
$A$} if $xg>x$ ($xg<x$) for some (and hence all) $x\in A$.
 
In \cite{BSPLR} and \cite{bpgsc} there are versions of the following
lemma and its corollary, but they can also be taken as exercises for
the reader.

\bl
\label{transitiveElementOrbital}

Suppose $g\in\ploi$ and $A = (a,\,b)$ is an orbital of $g$, then given
any $\epsilon>0$ and $x\in(a,\, b)$ there is an integer $n\in\Z$ so that
the following two statements are true.  

\be
\item $xg^n - a < \epsilon$, and
\item $b - xg^{-n} < \epsilon$.
\ee 
\el 

The following corollary of the above lemma can be proved by using a
compactness argument.

\bc
\label{transitiveOrbital}

Suppose $G\leq\ploi$ and $G$ has an orbital $A$. If $x<y$ are in $A$
then there is an element $\theta\in G$ so that $x\theta >y$.  
\ec

If $g\in \ploi$, and $A = (a,b)$ is an open interval in $I$, then we
will say that \emph{$g$ has an orbital that shares an end with $A$} if
$g$ has an orbital of the form $(a,c)$ or $(c,b)$.  We will also say
that \emph{$g$ realizes an end of $A$} in these situations.

The main result of \cite{BrinU} is given below.

\bt [Ubiquitous F]
\label{ubiquitousF}

If $H\leq\ploi$, and $H$ has an orbital $A$ so that some element $h\in
H$ realizes one end of $A$, but not the other, then $H$ contains a
subgroup isomorphic to Thompson's group $F$.  
\et

We will say a subgroup $G$ of
$\ploi$ is \emph{balanced} (following the language in \cite{bpasc}) if
it has no subgroup $H$ which satisfies the hypotheses of Theorem
\ref{ubiquitousF}.  We will see some properties of balanced groups
further on in this section.

For the statement of the next lemma to make sense, we need more
definitions.  If $g\in\ploi$, we will define $\brkp{g}$, the set of
breakpoints of $g$, to be the set of points in $(0,1)$ where the
derivative of $g$ is not defined.  We will call the components of
$[0,1]\backslash\brkp{g}$ the \emph{affine components of $g$} (note that
these are simply subsets of the domain of $g$ over which $g$ is
affine).  If $A$ is a connected subset of $[0,1]$, and $C$ is an
affine component of $g$, then we will call $C\cap A$ an \emph{affine
component of $g$ in $A$}.  The following is a restatement of Remark
2.1 in \cite{bpgsc}.

\bl 
\label{conjProps}

Suppose $g$,
$h\in\ploi$and that $A$ is an orbital of $g$.  The following are true.
\be
\item $g^h$ has orbital $Ah = \left\{ah\,|\, a\in A\right\}$.
\item If $g$ moves points right (resp. left) on $A$, then $g^h$ moves
points right (left) on $Ah$.
\item The slope of $g$ on the rightmost (resp. leftmost) affine
component of $g$ in $A$ equals the slope of $g^h$ on the rightmost
(leftmost) affine component of $g^h$ in $Ah$.
\ee
\el

Note that it is a straightforward consequence of the previous lemma
that the orbitals of $g$ and $g^h$ are in one-to-one ordered (left to
right) correspondence.  In general, we will refer to the orbital $Ah$
as the \emph{orbital of $g^h$ induced from $g$ by the action of $h$}
(or by other similar language).

Suppose $G$ is a subgroup of $\ploi$ and $A = (a,\,b)$ is an orbital
of $G$, further suppose there is an element $g$ that realizes both
ends of $A$.  Suppose the slope of $g$ on the leftmost affine
component of $g$ which non-trivially intersects $A$ is $s_l$, while
the slope of $g$ on the rightmost affine component of $g$ which
non-trivially intersects $A$ is $s_r$.  We say that \emph{$g$ realizes
$A$ inconsistently} if $s_l$ and $s_r$ are either both greater than
one or both less than one.  Otherwise, we say that \emph{$g$ realizes
$A$ consistently}.  If $g$ has $A$ as an orbital, then we say that $g$
realizes $A$ (note that in this case, $g$ realizes $A$ consistently).

We will now mention some promised properties of balanced subgroups of $\ploi$.
The interested reader is encouraged to examine section 3.3 of
\cite{bpasc}.  The following lemma lists two straightforward
consequences of the discussion there.

\bl
\label{consistentRealization}
\label{inconsistentRealization}

Suppose $G$ is a balanced subgroup of $\ploi$ and $G$ has an orbital
$A$ and an element $g$ which realizes both ends of $A$.
\be
\item  If $g$ realizes both ends of $A$ consistently, then $g$ realizes $A$.
\item If $g$ realizes $A$ inconsistently, then no element of $\,G$
realizes $A$ consistently.
\ee 

\el

Given a tower $T$ for a group $G\leq\ploi$, we may pass to the group
$\langle S_T\rangle$ generated by the signatures of the tower.  This
group can be fairly complicated, depending in part on how the other
orbitals of the signatures of $T$ align with each other.  We will say
a tower $T$ is an \emph{exemplary tower} if whenever $(A,g)$,
$(B,h)\in T$ with $A\neq B$ then $(A,g)\leq (B,h)$ implies both

\be
\item the orbitals of $g$ are disjoint from the ends of the
orbital $B$, and
\item no orbital of $g$ in $B$ shares an end with $B$.  
\ee

The following is a conglomeration of results from \cite{bpgsc} (Remark
2.9, Lemma 2.12, Lemma 3.2, and Remark 4.1) with a straightforward
extension of Lemma 3.14.1 in \cite{bpasc}.

\bl 
\label{productOrbitals} 

Suppose $G$ is a subgroup of $\ploi$ that fails to admit
transition chains of length two.  We have the following consequences.
\be
\item $G$ is balanced.
\item If $\,T$ is a tower for $G$, then $T$ is exemplary.
\item If $f$ and $g\in G$, and $A$ is an orbital of $f$ and $B$ is an
orbital of $g$, and $A\cap B\neq \emptyset$, then one of the following
three statements holds.  

\be
\item $A = B$, and $A$ is an orbital of $\langle f, g\rangle$.

\item $\overline{A}\subset B$, $A\cap Ag = \emptyset$, and $B$ is an
orbital of $fg$, $gf$, and of $\langle f,g\rangle$.

\item $\overline{B}\subset A$, $B\cap Bf = \emptyset$, and $A$ is an
orbital of $fg$, $gf$, and of $\langle f, g \rangle$.  

\ee 
\item If $G$ is a subgroup of depth $n$ for some positive integer $n$,
  then $G'$ has depth $n-1$.

\ee 
\el

We sometimes need to understand whether an element-orbital survives as
an element-orbital in a derived group.  

\bl 
\label{commutatorOrbital}

Suppose $G$ is a subgroup of $\ploi$.  If $T = \left\{(A_1,g_1),
(A_2,g_2)\right\}$ is an exemplary tower of height two for $G$, where
$(A_1,g_1)\leq (A_2,g_2)$, then there is $M\in\N$ so that for all
$n\in \N$ with $n \geq M$ we have that $A_1$ is an orbital of the
element $[g_1,g_2^n]$.  \el
pf:

Suppose $A_2 = (a,b)$, and let $[x,y]$ be the smallest interval so
that $\supp{g_1}\cap A_2 \subset [x,y]$. (This interval exists by the
definition of an exemplary tower.)  We see immediately that
$A_1\subset [x,y]$.  Let $\epsilon = \max(x-a,b-y)$.  By
\ref{transitiveElementOrbital} there is $S\in\Z$ so that $xg_2^S>y$.
Let $M = |S|$.  Note that for any $n>M$, we have $[x,y]g_2^n\cap[x,y]
= \emptyset$.  The lemma follows immediately.

\qquad$\diamond$   

  Suppose two elements $h$, $k\in\ploi$ have the property that
whenever $A$ is an orbital of $h$ and $B$ is an orbital of $k$ so that
$A\cap B\neq\emptyset$, then either $A = B$, $\overline{A}\subset B$,
or $\overline{B}\subset A$.  We say that $h$, $k$ are
\emph{mutually efficient}, or that they satisfy the \emph{mutual efficiency
condition}, if given any orbital $C$ of $h$ that contains the closure
of an orbital of $k$, then the support of $k$ in $C$ is contained in a
single fundamental domain of $h$ in $C$, and the symmetric condition
that whenever $D$ is an orbital of $k$ that properly contains the
closure of an orbital of $h$, then the support of $h$ in $D$ is
contained in a single fundamental domain of $k$ in $D$.

Given two mutually efficient elements $h$ and $k$ in $\ploi$, we often
will form the commutator $[[h,k],k]$ (recall that we use the
definition $[a,b] = a^{-1}b^{-1}ab$ for the commutator symbol), which
we will refer to as the \emph{double commutator of $h$ and $k$}.  The
following is a restatement of Lemma 4.2 in \cite{bpgsc}.

\bl
\label{dcFacts}

Let $h$, $k\in H$, where $H$ is a subgroup of $\ploi$ with no
transition chains of length two.  Suppose further that $h$ and $k$ are
mutually efficient.  If $f = [[h, k],k]$, then $f$ has the following
properties:

\be

\item Every orbital of $h$ whose closure is contained in an orbital of $k$ is
an orbital of $f$.

\item Every orbital of $f$ has closure contained in an orbital of $k$ that
contains (perhaps not properly) an orbital of $h$.

\ee 

\el

From this point forward, all results and discussion will be new.

\section{Relationships Amongst the Key Groups \label{wreathProducts}}

Let us now investigate the relationships between our key groups.  In
the discussions below, we will call $\beta_{-1}$ the top generator of
$^{\infty}(\wr\Z)$, and $\beta_0$ the bottom generator of
$(\Z\wr)^{\infty}$.

\bl 
\label{WInStuff}

$W$ embeds in each group in the set $\{^{\infty}(\wr\Z),\, (\Z\wr)^{\infty},\, (\wr\Z\wr)^{\infty},\,B\}$.  
\el

pf: 

We will show that $W$ embeds in $(\Z\wr)^{\infty}$ and
$^{\infty}(\wr\Z)$.  This will complete the proof since
$^{\infty}(\wr\Z)$ and $(\Z\wr)^{\infty}$ each embed in
$(\wr\Z\wr)^{\infty}$ and $B$.

We first embed $W$ in $(\Z\wr)^{\infty}$.  For each $j\in\N$, define  
\[
\Gamma_j = \left\{\gamma_{i,j}\,|\, \gamma_{i,j} = \beta_i^{\beta_{j+1}},
1\leq i\leq j, i\in\N\right\}
\]
Note that each collection $\Gamma_i$ generates a group isomorphic to
$W_i$, by the argument given in the introduction after the discussion
of $W_2$, or also by the details of Brin in \cite{BrinEG} (re-define
$\Gamma_0$ to be the set with only the identity element of $\ploi$).
Further, the supports of the generators in $\Gamma_i$ are all disjoint
from the supports of the generators in $\Gamma_j$ whenever $i$,
$j\in\N$ with $i\neq j$, so that the elements of $\Gamma_i$ found in
$(\Z\wr)^{\infty}$ above commute with the elements of $\Gamma_j$ in
this case.  Hence, the set
\[
\Gamma = \cup_{i\in\N}\Gamma_i
\]
generates a group
\[
\langle \Gamma \rangle\cong \bigoplus_{i\in\N}W_i \cong W.\]

(Note: We will use this realization of $W$ throughout the rest of this
subsection when we refer to our realization of $W$ in $\ploi$.  When
we refer to ``the first $n$ summands of $W$'' we will mean the
subgroup $\langle \cup_{i= 0}^n\Gamma_i\rangle\cong\bigoplus_{i = 0}^n
W_i$ of $W$. (Note that we are ignoring the trivial $W_0$
summand in our count.)

We now embed $W$ in $^{\infty}(\wr\Z)$ in a similar fashion, finding
copies of each $W_i$ in $^{\infty}(\wr\Z)$, all of which occur with
mutually disjoint supports in $I$, the union of their generators will
then generate a group isomorphic to $W$. Let $i\in\N$, and define
\[
\Upsilon_i=\left\{\theta_{i,j}\,|\, \theta_{i,j} =
\beta_{-i+j-2}^{\beta_{-1}^i}, 1\leq j\leq i, j\in\N\right\}.
\]
so that $\Upsilon_i$ is the collection of the $i$'th conjugates of the
$i$ generators beneath $\beta_{-1}$ of the generators of
$^{\infty}(\wr\Z)$.  Each collection $\Upsilon_i$ therefore generates
a group isomorphic to $W_i$ (re-define $\Upsilon_0$ to be
the set containing only the identity element of $\ploi$), while if
$i\neq j$, any generator in $\Upsilon_i$ has disjoint support from the
generators of $\Upsilon_j$, so that the union
\[
\Upsilon =
\cup_{i\in\N}\Upsilon_i
\]
 has the property that
\[
\langle \Upsilon\rangle\cong \bigoplus_{i\in\N}W_i\cong W.
\]
\qquad$\diamond$

Since $W$ can be realized as a subgroup of $\ploi$, the following
lemma demonstrates that the answer to Brin's Question 1 in
\cite{BrinEG} is ``No.''  

\bl
\label{stuffInW} 

No group in the set $\{^{\infty}(\wr\Z),\, (\Z\wr)^{\infty},\,
(\wr\Z\wr)^{\infty},\,B\}$ embeds in $W$.

\el 

pf:

We show that neither $^{\infty}(\wr\Z)$ nor $(\Z\wr)^{\infty}$ embeds
in $W$.  This will imply that $(\wr\Z\wr)^{\infty}$ and $B$ both fail
to embed in $W$.

There is a short proof of this restricted statement using the work of
Brin in \cite{BrinEG}, and our previous lemma.  By our previous lemma,
if $(\wr\Z)^{\infty}$ embeds in $W$, then it must embed in
$(\Z\wr)^{\infty}$.  Likewise, if $(\Z\wr)^{\infty}$ embeds in $W$,
then it must embed in $(\wr\Z)^{\infty}$.  However, neither of these
two things can occur; Theorem 2 in \cite{BrinEG} states that the two
groups $(\wr\Z)^{\infty}$ and $(\Z\wr)^{\infty}$ do not embed in each
other.
\qquad$\diamond$

\subsection{Transition Chains and the group $B$}

In this section we show that $B$ embeds in any subgroup of $\ploi$
which admits transition chains of length two.

\btl
\label{spanningConjugate} 

Suppose $H$ is a balanced subgroup of $\ploi$ and $H = \langle\alpha,
\beta\rangle$ for some two elements $\alpha$, $\beta\in\ploi$.
Suppose further that $A$ is an inconsistent orbital of $H$ and
$\alpha$ realizes both ends of $A$ while $\beta$ realizes neither.
There is a conjugate $\gamma$ of $\beta$ in $H$ which has an orbital
$B\subset A$ so that the fixed set of $\alpha$ in $A$ is contained in
$B$.

\etl

Pf: Let $F_{\alpha}$ represent the fixed set of $\alpha$ in $A$, and
let $x = \inf(F_{\alpha})$ and $y = \sup(F_{\alpha})$.  By Lemma
\ref{transitiveOrbital}, since $A$ is an orbital of $H$, there is
$\theta\in H$ so that $x\theta >y$.  By the continuity of $\theta$,
there is $x_1<x$ so that $x_1\theta>y$ as well.  Let $z = x_1\theta^{-1}$ so
we have
\[
z<z\theta = x_1<x<y<x_1\theta<x\theta
\]
Since $F_{\alpha}$ is contained in the orbitals of $\beta$, we see that $\beta$
has an orbital $C = (r,s)$ so that $r<x<s$.  By Lemma \ref{transitiveElementOrbital} there is a power $k\in\Z$
so that $r\alpha^k=q<z$.  Now, $\beta_1 = \beta^{\alpha^k}$ has
orbital $D=(q,t)$ induced from $C$ by the action of $\alpha^k$, and
$D$ satisfies $q<z<x<t$.  Set $\gamma = \beta_1^{\theta}$.  $\gamma$
has orbital $B =(u,v)$ induced from $D$ by the action of $\theta$ on
$\beta_1$, and $u = q\theta<x<y<t\theta = v$.

\qquad$\diamond$

The following lemma is more involved, and plays a key role in
the proof of the lemma following immediately after.  

\btl
\label{chainSplitting}

Suppose $H$ is a balanced subgroup of $\ploi$ and $H = \langle
\alpha,\beta\rangle$ for some two elements $\alpha$, $\beta\in\ploi$.
If $H$ has an inconsistent orbital $A$, and $\beta$ realizes no end
of any orbital of $H$, then there are elements $\alpha_1$ and
$\beta_1$ in $H$ so that if $H_1 = \langle\alpha_1,\beta_1\rangle$ the
following will be true.

\be

\item $A$ is an orbital of $H_1$.
\item $\beta_1$ realizes no end of any orbital of $H_1$.
\item Every inconsistent orbital of $H_1$ can be
written as the union of the orbitals of a transition chain of length
three, whose first and last orbitals are orbitals of $\alpha_1$, and
whose second orbital is an orbital of $\beta_1$.
\item $\alpha_1$ moves points to the left on its leading orbital in
each of the inconsistent orbitals of $H_1$.

\ee
\etl

pf: 
We break the proof into stages, so as to make it less cumbersome.

\vspace{.1 in}
\underline{S1: Classifying orbital types}.

Set $\alpha_1$ either to $\alpha$ or $\alpha^{-1}$, so that
$\alpha_1$ moves points to the left on its leading orbital $B$
contained in $A$.

Suppose $n\in\N$ and $H$ has $n$ inconsistent orbitals.  Let
$\mathscr{B}=\left\{B_i\,|\,1\leq i\leq r\right\}$ represent the
collection of inconsistent orbitals of $H$ where $\alpha_1$ moves
points to the left on its leading orbital in each of these orbitals,
indexed from left to right, where $r$ is the total number of such
orbitals.  Let $\mathscr{C}=\left\{C_j\,|\,1\leq j\leq s\right\}$
represent the other inconsistent orbitals of $H$, so that $n = s+r$,
where these orbitals are indexed from the left to the right as before.
Note that it is possible for $s=0$, although of course $r\geq 1$.

\vspace{.1 in}

\underline{S2: Building an element to span the fixed sets of $\alpha_1$ in the $B_i$}.

By Technical Lemma \ref{spanningConjugate}, for each orbital $B_i$ in
$\mathscr{B}$ there is an element $\gamma_i$ in $H$, which is a
conjugate of $\beta$, so that the fixed set of $\alpha_1$ in $B_i$ is
contained in a single orbital of $\gamma_i$.  Likewise, for each
orbital $C_j$ in $\mathscr{C}$ there is an element $\theta_j$ in $H$,
which is a conjugate of $\beta$, so that the fixed set of $\alpha_1$
in $C_j$ is contained in a single orbital of $\theta_j$.

Firstly, inductively replace each element $\gamma_i$, for $i>1$, by a
conjugate of $\gamma_i$ by a high negative power of $\alpha_1$ so that
for each $j$ with $1\leq j < i$ the closure of the union of all of the
orbitals of $\gamma_i$ in $B_j$ is actually fully contained in the
single orbital of $\gamma_j$ that contains the fixed set of $\alpha_1$
in $B_j$.  We can do this due to the specified directions in which 
$\alpha_1$ moves points on its first and last orbitals in each of the
$B_k$.

Summing up, for each $i\in\left\{1,2,\ldots,r\right\}$, $\gamma_i$ has
an orbital $D_i$ that contains the fixed set of $\alpha_1$ in $B_i$,
as well as the closure of all of the orbitals of $\gamma_k$ in $B_i$
for all $k\in \N$ where the inequalities $i<k\leq r$ hold.

We will now inductively define a sequence of elements $(\rho_i)_{i =
1}^r$ so that the $\rho_i$ will have the following properties (modulo
the fact that some of the $\gamma_i$ below will actually be conjugates
of the existing $\gamma_i$ by further negative powers of $\alpha_1$):

\be

\item $\rho_1 =\gamma_1$.
\item For all indices $i>1$, $\rho_i$ will be either a conjugate
of $\rho_{i-1}$ by some power of $\gamma_i$, or
$\rho_i=\rho_{i-1}\gamma_i$.
\item For all indices $i$, $\rho_i$ will have an orbital $E_i$ in
$B_i$ that fully contains the fixed set of $\alpha_1$ in $B_i$.
\item If $i<r$, the orbital $E_i$ of $\rho_i$ will contain the
  closures of the orbitals of $\gamma_j$ in $B_i$ for all integers 
  $j$ with $i < j \leq r$.
\item If $i>1$, for each integer $j$ with $1\leq j <i$, $\rho_i$ will
have $E_j$ as one of its orbitals.

\ee 

Firstly, set $\rho_1 =\gamma_1$, and $E_1 = D_1$.  By construction,
$\rho_1$ satisfies the five inductive properties.  If $r = 1$, we are
done.  If not, suppose that $k$ is an integer so that $1<k\leq r$ and
for all $i\in\left\{1,2,\ldots,k-1\right\}$ we have that $\rho_i$ is
defined and satisfies the five defining properties of the induction.
Our analysis now breaks into two cases.

If $\rho_{k-1}$ has an orbital $F_k$ containing either end of $D_k$,
then there is some integer $j$ so that $\rho_k =
\rho_{k-1}^{\gamma_k^j}$ will have orbital $E_k$ induced from $F_k$ by
the action of $\gamma_k^j$ so that $E_k$ will contain the fixed set of
$\alpha_1$ in $B_k$, as well as the closure of all of the orbitals of
$\gamma_j$ in $B_k$ for integers $j$ where the inequalities $k<j\leq
r$ hold.

If $\rho_{k-1}$ does not have an orbital $F_k$ containing either end
of $D_k$, then we have to handle the case where $\rho_{k-1}$ has
orbitals in $D_k$ that share ends with $D_k$ separately before
continuing.

If $\rho_{k-1}$ has orbitals in $D_k$ that share ends with $D_k$ then
replace $\gamma_k$ and all later $\gamma_j$ with conjugates of these
elements by a high negative power of $\alpha_1$ so that $D_k$ either
has an end contained in an orbital of $\rho_{k-1}$, or shares no end
with an orbital of $\rho_{k-1}$, and repeat the whole inductive
definition of $\rho_k$.

If $\rho_k$ is still undefined, then set $\rho_k =
\rho_{k-1}\gamma_k$.  Note that since $\rho_{k-1}$ has no orbitals in
$D_k$ that share ends with $D_k$, the product $\rho_k =
\rho_{k-1}\gamma_k$ realizes both ends of $D_k$ consistently, and
therefore realizes $D_k$ consistently since $H$ is balanced.
Therefore define $E_k = D_k$ and note that $\rho_k$ actually has $E_k$
as an orbital.

At this point, $\rho_k$ and $E_k$ are both defined, and we can
continue with our main argument.  Note that $E_k$ contains the closure
of all of the orbitals of all of the $\gamma_i$ for $i>k$, and that
for each integer $j$ in $1\leq j\leq {k-1}$, the closure of the
orbitals of $\gamma_k$ in $B_j$ are fully contained in the orbital
$E_j$, so that by Lemma \ref{consistentRealization} $\rho_k$ will have
$E_j$ as an orbital as well.  Now by construction, $\rho_k$ satisfies
the five defining properties of the induction.

We now examine the element $\rho_r$.  Observe that the element
$\rho_r$ contains an orbital $E_k$ in each $B_k$ where the fixed set
of $\alpha_1$ in $B_k$ is fully contained in $E_k$.  $\rho_r$ is
constructed as a sequence of products using various $\gamma_i$'s and
conjugates of $\gamma_i$'s so $\rho_r$ realizes no end of any orbital
of $H$, but is an element of $H$.

\vspace{.1 in}

\underline{S3: Building an element to span the fixed sets of
$\alpha_1$ in the $C_j$}.

In an entirely analogous fashion, if $s>0$, then we can find one
element $\psi_s$ in $H$ which realizes no end of any orbital of $H$
and which contains an orbital $F_i$ in each $C_i\in\mathscr{C}$ which
contains the fixed set of $\alpha_1$ in that $C_i$.

\vspace{.1 in}

\underline{S4: Modifying our first element so that it creates no transition
chains with $\alpha_1$ over the $C_i$}.

There is a positive integer $p$ so that $\rho= \rho_r^{\alpha_1^p}$
has the following two properties.

\be

\item For each integer $i\in\left\{1,2,\ldots,r\right\}$, the closure
of the orbitals of $\psi_s$ in $B_i$ is actually contained in the
orbital $G_i$ of $\rho$ induced from $E_i$ by the action of
$\alpha_1^p$.
\item For each integer $i\in\left\{1,2,\ldots,s\right\}$, the closure
of the orbitals of $\rho$ in $C_i$ is actually contained in the
orbital $F_i$ of $\psi_s$.

\ee

This follows since for each orbital $B_i$ of $\mathscr{B}$, the lead
orbital of $\alpha_1$ in $B_i$ has the property that $\alpha_1$ is
moves points to the left there (and therefore moves points to the
right on the trailing orbital of $\alpha_1$ in $B_i$), and for each
orbital $C_i$ in $\mathscr{C}$, the lead orbital of $\alpha_1$ in
$C_i$ has the property that $\alpha_1$ is moving points to the right
there (and therefore $\alpha_1$ moves points to the left on its
trailing orbital in $C_i$).

We note in passing that the orbitals $G_i$ of $\rho$ contain the
orbitals $E_i$ of $\rho_r$, and therefore the fixed set of $\alpha_1$ in
the $B_i$.

Now there is a power $q$ of $\psi_s$ so that the element $\beta_1 =
\rho^{\psi_s^q}$ will have the following nice properties: 

\be

\item For each integer $i\in\left\{1,2,\ldots,s\right\}$, the orbitals
of $\beta_0$ in $C_i$ have trivial intersection with the fixed set of
$\alpha_1$ in $C_i$.

\item For each integer $i\in\left\{1,2\ldots,r\right\}$, $\beta_0$ will
have the orbital $G_i$ which contains the fixed set of $\alpha_1$ in
$B_i$.

\ee 

The first property follows since the orbitals of $\rho$ in the $C_i$
are contained in the orbitals $F_i$ of $\psi_s$, and so the
conjugation of $\rho$ by a high power of $\psi_s$ will throw these
orbitals off of the fixed set of $\alpha_1$ in the $C_i$.  The second
property follows since the orbitals of $\psi_s$ are fully contained in
the orbitals $G_i$ of $\rho$ in the $B_i$, so that conjugation of
$\rho$ by $\psi_s$ to any power will not change these orbitals.

It is now straightforward to check that the group $H_1 = \langle
\alpha_1,\beta_1\rangle$ satisfies all of the properties promised in the
statement of the lemma.
\qquad$\diamond$ 

We are now ready to prove our chief geometric result.  Below is a
re-statement of Theorem \ref{tChainB}.

\bt

If $G$ admits a transition chain of length two, then $B$ embeds in $G$.

\et

pf:

We can assume that $G$ is balanced, otherwise $G$ will contain a copy
of R. Thompson's group $F$, which contains copies of Brin's group $B$.

Let $\mathscr{T} = \left\{(O_1,\alpha), (O_2,\beta)\right\}$ be a
transition chain of length two for $G$ and let $K = \langle
\alpha,\beta\rangle$.

The orbitals of $K$ are the components of the union of the orbitals of
$\alpha$ and $\beta$.  Some of these orbitals may be consistent
orbitals for $K$, so that at least one of $\alpha$ or $\beta$ realize
these orbitals.  The other orbitals are inconsistent, and are formed
by the union of a sub-collection of orbitals of $\alpha$ and orbitals
of $\beta$.  A chief feature of the inconsistent orbitals is that one
of $\alpha$ or $\beta$ must realize both ends of any particular such
orbital, since $K$ is balanced.  Note in passing that the orbital of
$K$ which contains $O_1$ is inconsistent, as neither $\alpha$ nor
$\beta$ realize that whole orbital, although one of $\alpha$ or $\beta$
realizes both ends of it.

We will proceed through the remainder of the proof in stages.

\vspace{.1 in}

\underline{S1: Classifying the orbitals of $K$}.

 We are going to analyze the orbitals of $K$ still further.  Any
 particular orbital of $K$ has one of six types, the first three are
 consistent, and the last three are inconsistent:

\be

\item (Type AB) Both $\alpha$ and $\beta$ consistently realize this orbital.
\item (Type Ab) $\alpha$ consistently realizes this orbital, but not $\beta$.
\item (Type aB) $\beta$ consistently realizes this orbital, but not $\alpha$.
\item (Type $\underline{a}\underline{b}$) Both $\alpha$ and $\beta$
inconsistently realize both ends of this orbital.
\item (Type $\underline{a}$b) $\alpha$ inconsistently realizes both
ends of this orbital, but $\beta$ realizes
neither end of this orbital.
\item (Type a$\underline{b}$) $\beta$ inconsistently realizes both
ends of this orbital, but $\alpha$ realizes
neither end of this orbital.

\ee

We know that $K$ has at least one orbital, let us call it $A$, of type
$\underline{a}\underline{b}$, $\underline{a}b$, or $a\underline{b}$,
and we will assume without meaningful loss of generality that $A$ has
one of the first two types.  Let $F_a$ represent the union of the
fixed sets of $\alpha$ that are contained in the orbitals of $K$ of
type $\underline{a}\underline{b}$, type $\underline{a}b$, and
$a\underline{b}$.  $F_a$ is non-empty, and is entirely contained in
the orbitals of $\beta$.  

\vspace{.1 in}

\underline{S2: A mechanism for modifying elements and orbitals}.

By Remark
\ref{transitiveElementOrbital} there is a $N_1\in\N$ so that for all
$k\in\N$ with $k\geq N_1$ we have $F_a\beta^k\cap F_a=\emptyset$ in
orbitals of type $\underline{a}\underline{b}$ and $\underline{a}b$
(also assume that $N_1$ is large enough so that in orbitals of type
$a\underline{b}$, the interior components of $F_a$ are moved off of
themselves of the action of $\beta^k$).  Similarly, let $S$ represent
the support of $\alpha$ in the orbitals of $K$ of type $aB$, then
there is $N_2\in\N$ so that for all $k\geq N_2$, we have $S\beta^k\cap
S=\emptyset$.  Let $N = max(N_1,N_2)$.  

Considering the other direction, let $F_b$ represent the fixed set of
$\beta$ in the orbitals of $K$ of type $a\underline{b}$.  Since $F_b$
is contained in the support of $\alpha$ by definition, there is
$M\in\N$ so that for all $j\geq M$ we have that $F_b\alpha^j\cap
F_b=\emptyset$.  

Now let $j\geq M$, and let $k\geq N$, and define
$\beta_1=[\alpha^j,\beta^k]$.  

\vspace{.1 in}

\underline{S3: Analyzing how the mechanism effected $\beta$}.

Observe that the fixed set of $\alpha$ in the orbitals of $K$ of type
$\underline{a}\underline{b}$ and $\underline{a}b$ is contained in the
orbitals of $\beta_1$.  Note also that the components of $F_a$ in the
orbitals of $K$ of type $a\underline{b}$ which do not realize any end
of an orbital of $K$ are all contained in the support of $\beta_1$.
The fixed set of $\beta$ contained in the orbitals of $K$ of type
$a\underline{b}$ is also contained in the support of $\beta_1$, since
any such point is moved off of $F_b$ by $\alpha^{-j}$, then moved by
$\beta^{-k}$, then moved to someplace different (from its start) by
$\alpha^j$, and finally, $\beta^k$ cannot move the resultant point to
its original location in the fixed set of $\beta$.  Now observe that
the orbitals of $\beta_1$ are either disjoint from $S$, or else are
components of $S$ where $\alpha^j$ behaves as the inverse of
$\beta_1$.

\vspace{.1 in}

\underline{S4: Analyzing the orbitals of $K_1$}.

We now consider the group $K_1=\langle \alpha, \beta_1\rangle$, and we
consider the orbitals of $K_1$ under the same classification as
the orbitals of $K$, where we replace $\beta$ by $\beta_1$ in that
classification.

It is immediate to see that $K_1$ still has all the orbitals of $K$ of
type $\underline{a}b$, and that the type of these orbitals is
unchanged.  It is also immediate by construction that the orbitals of
$K$ of type $\underline{a}\underline{b}$ are also orbitals of $K_1$,
although they are now of type $\underline{a}b$.  The orbitals of $K$
of type $Ab$ are also orbitals of $K_1$ of type $Ab$, but the orbitals
of $K$ of type $aB$ are now replaced by a collection of interior
orbitals (all lying properly in the union of the orbitals of $K$ of
type $aB$), each of which is an orbital of type $aB$ that is actually
disjoint from the support of $\alpha$, or else of type $AB$, where
$\alpha^j$ and $\beta_1$ behave as inverses on these orbitals.  The
orbitals of $K$ of type $AB$ are now of type $Ab$, and may have
trivial intersection with the support of $\beta_1$ (if, in fact,
$\alpha$ and $\beta$ commuted on these orbitals).

If $B_1$ is an orbital of $K$ of type $a\underline{b}$, then $B_1$ is
not an orbital of $K_1$.  In this case $K_1$ admits a new collection
of orbitals properly contained in $B_1$.  

We first consider the case where $\beta$ is moving points to the right
on its leading orbital in $B_1$ (and therefore is moving points to the
left on its trailing such orbital). We will suppose $k$ was chosen
large enough so that the closure of the union of the orbitals of
$\beta^{-k}\alpha^j\beta^k$ that are contained in orbitals of $\beta$
in $B_1$ is actually contained in the orbitals of $\alpha$ (and
therefore of $\alpha^j$) which contain components of the fixed set of
$\beta$.  Note that any interior orbital of $\beta$ in $B_1$ is
contained in the union of the orbitals of $\alpha^j$ and
$\beta^{-k}\alpha^j\beta^k$.  Therefore, there are three possible
varieties of resulting orbitals of $K_1$ in $B_1$: firstly, of type
$AB$, where $\beta_1$ actually behaves as $\alpha^{-j}$ on these
orbitals (there may be several of these), secondly, of type $AB$,
where there is only one such orbital, and it contains the fixed set of
$\beta$, or thirdly, of type $\underline{a}\underline{b}$, where there
is one of these if the previous variety did not occur, and it contains the
fixed set of $\beta$ in this case.  We will assume $k$ was chosen
large enough so that these properties of transformation are preserved
over all orbitals of $K$ of type $a\underline{b}$ where $\beta$ is
moving points to the right on its leading relevant orbitals.

In the case of the orbitals of $K$ of type $a\underline{b}$ where
$\beta$ moves points to the left on its leading relevant orbitals. The
results depend heavily on the nature of $\alpha$ in these individual
orbitals.  To clarify the discussion, let us suppose that $B$ is such
an orbital, and discuss the possibilities that arise from the
behavior of $\alpha$ and $\beta$ on $B$.  

Firstly, let us suppose that $\alpha$ has an orbital that contains the
fixed set of $\beta$ in $B$.  In this case, let us suppose $k$ and $j$
were chosen large enough so that the entire support of $\alpha$ is
contained inside a single fundamental domain of the single orbital of
$\beta^{-k}\alpha^j\beta^k$ that contains the fixed set of $\beta$ in
$B$.  In this case, the group $K_1$ possibly has several orbitals in
$B$, all of type $aB$.  One of these orbitals contains all of the
support of $\alpha$ in $B$, and all of the rest are orbitals of
$\beta_1$ which contain no orbitals of $\alpha$ and are therefore of
type $aB$ with trivial intersection with orbitals of $\alpha$.

Now let us suppose that $\alpha$ has more than one orbital in $B$ that
contains a component of the fixed set of $\beta$.  The first and last
such orbitals of $\alpha$ in $B$ must have that $\alpha$ behaves
inconsistently on these orbitals, otherwise it is easy to create an
imbalanced subgroup of $K_1$.  So now there are two further cases.  

Let us suppose that $\alpha$ moves points to the right on its first
orbital in $B$ which contains a component of the fixed set of $\beta$,
and therefore moves points to the left on the last such.  In this case
$K_1$ has only one orbital in the domain $B$, call it $C$, which is
again of type $a\underline{b}$.  The closure of $C$ is contained in
$B$, and $\beta_1$ moves points to the right on its leading orbital in
$C$ and moves points to the left on its trailing orbital there.  

Now let us suppose $\alpha$ moves points to the left on its leading
orbital in $B$ that contains a component of the fixed set of $\beta$,
and therefore moves points to the right on its trailing orbital in $B$
which contains a component of the fixed set of $\beta$, the group
$K_1$ again has some pure orbitals (type $aB$) plus precisely one
orbital $C$ in $B$, which is again of type $a\underline{b}$, and this
time, $\beta_1$ will move points to the left on its leading orbital in
$C$ and will move points to the right on its trailing orbital in $C$.

The result of all of this analysis is the following, we can choose $j$
and $k$ so that the group $K_1$ has orbitals of the following types:

\be
\item $AB$

Note that in this case $\alpha$ and $\beta_1$ commute on this orbital,
except in the case possibly generated from orbitals of type
$a\underline{b}$ where $b$ moves points right on its leading orbital.

\item $Ab$
\item $aB$

Note here that the behavior of $\alpha$ on this orbital is as the
identity, unless this orbital is contained in an orbital of $K$ of
type $a\underline{b}$, in which case $\alpha$ may have non-trivial
support in this orbital.

\item $\underline{a}\underline{b}$ 

Note that orbitals of this type are always contained in orbitals of
$K$ of type $a\underline{b}$ where $\beta$ moves points to the
right on its first relevant orbital.

\item $\underline{a}b$

Since these are the certain result of an orbital of type
$\underline{a}\underline{b}$ or of an orbital of type $\underline{a}b$
of $K$, we see that $K_1$ will have at least one of these.

\item $a\underline{b}$

These orbitals all have the property that whenever $\beta_1$ moves
points to the left on its leading orbital in these orbitals, then
$\alpha$ moves points to the left on its leading orbital of the
orbitals that contain a component of the fixed set of $\beta_1$.

\ee

\vspace{.1 in}

\underline{S5: Evolving the group $K$ and its orbitals by repeatedly
applying the mechanism}.

We can repeat the process above to create a new element $\beta_2$
using $\alpha$ and $\beta_1$, and therefore a new group $K_2 =
\langle\alpha,\beta_2\rangle$.  $K_2$ improves on $K_1$ since all of
its orbitals of type $a\underline{b}$ have both $\beta_2$ and $\alpha$
moving points to the left on their important leading orbitals.  In
particular, $K_2$ may still have orbitals of type
$\underline{a}\underline{b}$, and of type $AB$ (although here $\alpha$
and $\beta_2$ will commute on these orbitals). $K_2$ may have orbitals
of type $Ab$, but its orbitals of type $aB$ will all have the property
that $\alpha$ is the identity over these orbitals, while $K_2$ will
certainly have orbitals of type $\underline{a}b$.  Repeating the
process one more time to create an element $\beta_3$ and a subgroup
$K_3=\langle \alpha,\beta_3\rangle$ produces a group whose orbitals
are much easier to describe.  $K_3$ will have no orbitals of type $AB$
since $K_2$ had no orbitals of type $aB$ or $a\underline{b}$ that
could produce these orbitals (the types exist, but not with the right
sub-flavors of $\alpha$ and $\beta_2$ to generate these offspring).
$K_3$ may have orbitals of type $Ab$, but it will have no orbitals of
type $aB$, since the orbitals of type $aB$ in $K_2$ had $\alpha$
behaving as the identity there, and $K_2$ had no orbitals of type
$a\underline{b}$ with $\beta_2$ moving points to the left on its first
sub-orbital $D$ while $\alpha$ was moving points to the right on its
orbital containing the right end of $D$.  $K_3$ will have no orbitals
of type $\underline{a}\underline{b}$, since $K_2$ had no orbitals of
type $a\underline{b}$ with $\beta_2$ moving points to the right on its
first orbital in the orbitals of $K_2$ of this type. $K_3$ will have
at least one orbital of type $\underline{a}b$, and may have several
orbitals of the type $a\underline{b}$, but all of these last will have
$\beta_3$ moving points to the left on its leftmost orbitals in these
orbitals, and $\alpha$ will also move points to the left on its first
orbitals containing the right ends of $\beta_3$'s leftmost orbitals in
these orbitals of type $a\underline{b}$ of $K_3$.

Now, the orbitals of $K_3$ are well understood, and the behaviors of
$\beta_3$ and $\alpha$ on these orbitals are also well understood.  We
now consider the subgroup $K_4$ generated by $\alpha$ and $\beta_4 =
[\alpha^{-j},\beta_3^{k}]$, where $j$ and $k$ are chosen as in the
previous process (note the negative index on $\alpha$).  The point of
this is that now the orbitals of $K_4$ will admit no orbital of type
$a\underline{b}$ with $\beta_4$ moving points to the left on its first
orbital.  Now replacing $K_4$ with $K_5 = \langle
\alpha,\,\beta_5\rangle$ where $\beta_5 = [\alpha^j,\beta^k]$ where
$j$ and $k$ are chosen as before produces a group with no orbitals of
type $a\underline{b}$, repeating one more time to generate $\beta_6$
and $K_6$ in the same fashion that we generated $K_1$ from $K$
produces a group whose orbitals are only of types $Ab$ and
$\underline{a}b$.

Let us consider the orbital $A$ of $K$.  $A$ is also an orbital of
$K_6$, and it is of type $\underline{a}b$. We will now replace $K$ by
$K_6$ and $\beta$ by $\beta_6$ so that $K$ has an orbital of type
$\underline{a}b$ and all of its orbitals are of type $\underline{a}b$
and $Ab$.

\vspace{.1 in}

\underline{S6: Improving the inconsistent orbitals of $K$}.

Suppose $K$ has $n$ orbitals of type $\underline{a}b$, and let
$\mathscr{O}=\left\{A_i\,|\,1\leq i\leq n, i\in\N\right\}$ represent
this collection, where the indices respect the left to right order of
the orbitals.  By construction we know that $n\geq 1$.  Apply
Technical Lemma \ref{chainSplitting} (above) to replace $\alpha$ and
$\beta$ by new elements, and replace $K$ by the new group generated by
the new $\alpha$ and $\beta$ so that $\beta$ still realizes no end of
any orbital of $K$, and $A_1$ is still an orbital of type
$\underline{a}b$, but where every maximal transition chain (of length
greater than one) which can be formed by using $\alpha$ and $\beta$
has length three (naturally $\alpha$ provides the leading and trailing
orbitals for any such chain), and where $\alpha$ moves points to the
left on all of its leading orbitals in orbitals of type
$\underline{a}b$ for $K$ (and therefore moving points to the right on
its trailing such intervals).

\vspace{.1 in}

\underline{S7: Improving the consistent orbitals of $K$}.

We now improve $\beta$ so that it will not admit support in the
orbitals of $K$ of type $Ab$.  

Choose two integers $m$ and $n$ intelligently.  Choose $n$ large
enough so that the entire support of $\beta$ in the inconsistent
orbitals of $K$ is contained in the set of orbitals of
$\beta^{\alpha^n}$ which contain the fixed sets of $\alpha$ in these
inconsistent orbitals, and so that in the consistent orbitals of $K$,
the conjugate $\beta^{\alpha^n}$ has support disjoint from the support
of $\beta$, being entirely to the left (or entirely to the right) of
the support of $\beta$ in each of these individual orbitals.  Choose
$m$ large enough so that the support of $\beta$ in the orbitals of
$\beta^{\alpha^n}$ (in the inconsistent orbitals of $K$) is moved
entirely to the right of itself (or entirely to the left of itself,
depending on the orbital of $\beta^{\alpha^n}$ involved) by the action
of $(\beta^m)^{\alpha^n}$.

Replace $\beta$ by the commutator $[\beta^m,(\beta^m)^{\alpha^n}]$.
The new beta still has a single orbital spanning the fixed set of
$\alpha$ in each of the inconsistent orbitals of the original $K$.
Also, the new $K$ generated by $\alpha$ and the new $\beta$ has the
same orbitals as the previous $K$, except now the support of $\beta$
is contained in the inconsistent orbitals of $K$.  In particular, the
consistent orbitals of $K$ now only support the action of powers of
$\alpha$.

\vspace{.1 in}

\underline{S8: Finding $B$ in $K$}.

Define $\tilde{\gamma}_0 = \beta$. There is a natural number $k$ so that
$\tilde{\gamma}_1 = \gamma_0^{\alpha^k}$ has the property that if $Z$ is any
particular orbital of orbital of $K$ of type $\underline{a}b$, then
the closure of the support of $\gamma_0$ in $Z$ is fully contained in
the orbital of $\tilde{\gamma}_1$ that contain the fixed set of $\alpha$ in
$Z$.  Note that this $k$ exists, since $\alpha$ moves points to the
left on all of its leading orbitals in the orbitals of $K$ of type
$\underline{a}b$.  Replace $\alpha$ by $\alpha^k$.

There is another natural number $j$ so that the closure of the support
of $\tilde{\gamma}_0$ in any particular orbital $Y$ of
$\tilde{\gamma_1}$ is fully contained in a single fundamental domain
of $\tilde{\gamma_1}^j$ in $Y$.  In particular, if we define $\gamma_0
= \gamma_0^j$, and $\gamma_i = \gamma_0^{\alpha^{i}}$ for every
$i\in\Z$, then the group generated by the $\gamma_i$ will be
isomorphic with $(\wr\Z\wr)^{\infty}$ since the support of each
$\gamma_i$ in any orbital of type $\underline{a}b$ of $K$ is wholly
contained inside a single fundamental domain of an orbital of
$\gamma_{i+1}$.  Since the orbitals of $\alpha_i$ which create the
consistent orbitals of $K$ do not effect the isomorphism type of the
group $\langle\alpha,\gamma_0\rangle$ (being disjoint from the
support of $\gamma_0$), we see that
$\langle\alpha,\gamma_0\rangle\cong B$.
 
\qquad$\diamond$

\subsection{Finding infinite wreath products in groups with infinite towers}

Suppose $D = \left\{(A_i,h_i)\,|\,i\in\N\right\}$ is an exemplary
tower whose indexing respects the order of the elements so that 
\be
\item $H = \langle S_D\rangle$ is a balanced group that admits no
transition chains of length two, and 
\item that whenever $B$ is an orbital of $h_i$ for some signature
$h_i$ of $D$, then $B$ is contained in an orbital $C$ of $h_{i+1}$.
\ee
We are going to find a sub-tower of $D$ that satisfies a nice further
property.

Suppose $B_1$ is an orbital of $h_1$.  Each signature $h_i$ of $D$ has
an orbital $B_i$ that contains $B_1$.  The orbitals $B_j$ are nested
as the index increases, but possibly not properly.  If there is an
$N\in\N$ so that for all $n>N$, we have $B_n = B_{n+1}$, then we
will call $B_1$ \emph{a terminal orbital of
$D$}\index{orbital!terminal}, and $(B_1,h_1)$ a terminal signed
orbital of $D$, and we will say that \emph{$B_1$ is stable after
$N_1$}.  We now extend this language to orbitals of signatures other than
$h_1$.  Given $i\in\N$, call an orbital of $h_i$ terminal in $D$ if
the orbital is terminal in the sub-tower of $D$ formed using only the
signed orbitals $(A_k,h_k)$ with $k\geq i$.  We will call any orbital
of a signature of the tower, where the orbital is not a terminal
orbital, a non-terminal orbital.  Observe that non-terminal signed
orbitals make good candidates for being bases of new exemplary towers.

We will rely heavily on the following technique in our proof of Lemma
\ref{tallZWreath}.

\bl[growing sub-tower] 
\label{growingTower}

Suppose $D = \left\{(A_i,h_i)\,|\,i\in\N\right\}$ is an exemplary
tower where $H = \langle S_D\rangle$ is a balanced group that admits
no transition chains of length two, and so that whenever $B$ is an
orbital of $h_i$ for some signature $h_i$ of $D$, then $B$ is
contained in an orbital $C$ of $h_{i+1}$.  Then we can pass to an
infinite sub-tower $E$ of $D$ so that if $J$ is an orbital of any
signature $g_i$ of $E$, where $J$ is not a terminal orbital of $D$,
then there is an orbital $K$ of $g_{i+1}$ which properly contains $J$.
\el

Pf:

We note by definition that the orbital $K$ will also be a non-terminal
orbital of $D$, and both will be non-terminal in $E$.

We now pass repeatedly to infinite sub-towers for $D$, at each stage
referring to the new tower that results as $D$, and re-indexing so
that the tower will still have the form $D =
\left\{(A_i,h_i)\,|\,i\in\N\right\}$.  Let $\mathscr{P} =
\left\{B_i\,|\,1\leq i\leq n_i, i\in\N\right\}$ represent the $n_i$
orbitals of $h_1$ that are not terminal, in left to right order.  We
improve $D$ by passing to a infinite sub-tower $n_i$ times.  Firstly,
for $B_1$, pass to a sub-tower of $D$ so that the orbitals of the $h_i$
over $B_1$ are always properly nested as we progress up the tower.
The new tower $D$ still has all the properties that we have listed for
the old $D$, but now the orbitals of the $h_i$ over $B_1$ actually
form a tower over $B_1$ when we pair them with their signatures.
Repeat this process inductively for each of the non-terminal orbitals
of $D$ in $h_1$.  Now we pass to an infinite induction, by repeating
the process again, using base signature $h_2$, so that we are
progressively improving the tower above $h_2$ so that the non-terminal
orbitals of $h_2$ are each actually the base of an infinite tower
using the signatures $h_k$ with $k>2$ paired with their appropriate
orbital containing the relevant orbital of $h_2$.  We note in passing
that the non-terminal orbitals of $h_1$ are all contained in the
non-terminal orbitals of $h_2$, so we only have to improve $D$ over
the non-terminal orbitals of $h_2$ which do not contain orbitals of
$h_1$.  With these observations in place, we can inductively continue
this process at every level of $D$.  Let $E$ be the tower that results
from this process.  Given any $i\in\N$, if $(A_i,h_i)$ is a signed orbital of $E$
and $B_i$ is a non-terminal orbital of $h_i$ then for any integer
$k>i$ there is an orbital $B_k$ of $h_k$ so that $\bar{B}_i\subset
B_k$.  Thus, $E_{B_i} = \left\{(B_k,h_k)\,|\,k\geq i, k\in\N\right\}$
is itself an exemplary tower.  \qquad$\diamond$

The following lemmas are simply restatements (with proofs) of the
different aspects of Theorem \ref{tallWreath}.

\bl
\label{tallZWreath}

If $G$ is a subgroup of $\ploi$ and $G$ admits a tall tower, then $G$
has a subgroup of the form $(\Z\wr)^{\infty}$.

\el

pf:

We again proceed in stages.

\vspace{.1 in}

\underline{S1: Making observations which enable a simplified treatment of orbitals}.

We will assume that $G$ is balanced, as otherwise $G$ contains a
subgroup isomorphic to Thompson's group $F$, which itself has a
subgroup isomorphic to $(\Z\wr)^{\infty}$.  If $G$ admits transition
chains of length two then by Theorem \ref{tChainB} $G$ admits an
embedded copy of Brin's group $B$, which also implies our result, so
let us assume that $G$ admits no transition chains of length two.

\vspace{.1 in}

\underline{S2: Choosing an initial tower, and noting its supporting orbital}.

Let $E = \left\{(A_i,g_i)\,|\,i\in\N\right\}$ be a tall tower for $G$,
where the indexing respects the order on the signed orbitals of $E$.
By Lemma \ref{productOrbitals}, since $G$ contains no transition
chains of length two, $E$ is exemplary.

Let $A = \cup_{i\in\N}A_i = (a,b)$.  We observe that if $B$ is an
orbital of $g_i$ for some $i$, then $B$ is disjoint from
$\left\{a,b\right\}$, and that if $B\cap A \neq\emptyset$, then
neither $a$ nor $b$ is an end of $B$.  In particular, $A$ is an
orbital of $\langle S_E \rangle$.

\vspace{.1 in}

\underline{S3: Improving our tower inside of the supporting orbital $A$}.

Now given $\epsilon>0$ so that $\epsilon< \frac{b-a}{2}$, we see that
there is an $N\in\N$ so that for all $n\in\N$ with $n\geq N$, we have
that $(a+\epsilon,b-\epsilon)\subset A_n$ since the ends of the $A_i$
must limit to the ends of $A$.  But now, we can construct a monotone
strictly increasing, order preserving function, $\phi:\N\to\N$, so
that given any $n\in\N$, all of the orbitals of $g_n$ in $A$ are
actually contained in $A_{\phi(n)}$, and since $E$ is exemplary, no
orbital of $g_n$ in $A$ actually shares an end with $A_{\phi(n)}$.  For any
$k\in\N$, let $\phi^k$ represent the product (via composition) of the
function $\phi$ with itself $k$ times in the monoid of order
preserving functions from $\N$ to $\N$ (use $\phi^0 = id$, the
function which moves nothing).  Now define an order
preserving function $\theta:\N\to\N$, defined by the rules that
$0\mapsto 0$ and $n\mapsto \phi^{n}(0)$ for each
$n\in\N\backslash\left\{0\right\}$.  Replace $E$ by the exemplary
tower formed by the collection
$\left\{(A_{\theta(i)},g_{\theta(i)})\,|\,i\in\N\right\}$.  

$E$ now has the property that if $i$, $k\in\N$ with $i<k$ then all the
orbitals of $g_i$ in $A$ are actually in $A_k$, away from the ends of
$A_k$.  For each $n\in\N$, with $n>0$, let $m_n$ be an integer large
enough so that the collection of orbitals of $g_{n-1}$ inside of $A_n$
(which is all the orbitals of $g_{n-1}$ in $A$) is actually fully
contained in a single fundamental domain of $g_n^{m_n}$ in $A_n$.
Define $m_0 = 1$.  Improve $E$ by replacing each signature $g_n$ with
$g_n^{m_n}$.  Now define $H = \left<S_E\right>$.  We note in passing
that $A$ is an orbital of $H$.

\vspace{.1 in}

\underline{S4: Improving our tower outside of the orbital $A$}.

We cannot immediately pass to a growing sub-tower; some work needs to
be done to $E$ in order for it to satisfy the hypotheses of Lemma
\ref{growingTower}.

Define $h_0=g_0$.  Now for each
$n\in\N$ with $n >0$, inductively define $h_n$ via the following four
step process.  

First, define $k_n = g_n^{r_n}$, where $r_n$ is a
positive integer large enough so that whenever $B$ is an orbital of
$g_n$ that is also an orbital of $h_{n-1}$, then the product
$h_{n-1}k_n$ still has orbital $B$.  

Second, define $h_n' = h_{n-1}k_n$.  Recall from Remark
\ref{productOrbitals} that any orbital of $h_{n-1}$ which properly
contains an orbital of $k_n$ will now be an orbital of $h_n'$,
and that any orbital of $k_n$ that properly contains an orbital of
$h_{n-1}$ will also be an orbital of $h_n'$.

$h_n'$ now has an orbital containing every orbital of $h_{n-1}$.

Third, choose positive integer $s_n$ large enough so that
every orbital $C$ of $h_n'^{s_n}$ which properly contains orbitals of
$h_{n-1}$ actually contains all such orbitals in a single fundamental
domain of $h_n'^{s_n}$ on $C$.  

Fourth, define $h_n = h_n'^{s_n}$.  The result is that the sequence
$(h_i)_{i\in\N}$ of signatures satisfies the following list of
properties.

\be

\item For each $n\in\N$, $A_n$ is an orbital of $h_n$.
\item For each $n\in\N$ with $n>0$, the orbitals of $h_{n-1}$ in $A$
are all contained inside a single fundamental domain of $h_n$ in
$A_n$.
\item For each $n\in\N$ with $n>0$, if $B$ is an orbital of $h_n$
which is not disjoint from the orbitals of $h_{n-1}$, then there are
two possibilities.

\be

\item $B$ is also an orbital of $h_{n-1}$.
\item $B$ properly contains a non-empty collection of orbitals of
$h_{n-1}$ in a single fundamental domain of $h_n$ on $B$.  

\ee 
\ee 

In particular, we can form the new exemplary tower $D =
\left\{(A_i,h_i)\,|\,i\in\N\right\}$.

$D$ still has the properties that $\cup_{i\in\N} A_i = A$, and that
$A$ is an orbital of the group $\langle S_D \rangle$.  Further, the
signatures satisfy the three enumerated points above.

We will now improve $D$ by replacing it with the result of finding
a growing sub-tower, so that any non-terminal orbital of any signature
$h_i$ of $D$ is properly contained in a non-terminal orbital of a
signature with index one higher.

\vspace{.1 in}

\underline{S5: Removing terminal orbitals}.

Our new $D$ is far superior to our old $D$, but $h_0$ will still have
terminal orbitals, if it had them to begin with.  Suppose $h_0$ does
have some terminal orbitals.  Then there is $N_0\in\N$ so that all the
terminal orbitals of $h_0$ are stable for $n\geq N_0$.  Compute a new
element $k = [[h_{N_0+1}, h_{N_0+2}], h_{N_0+2}]$ (note that condition
(3) above implies that $h_{N_0+1}$ and $h_{N_0+2}$
satisfy the mutual efficiency condition since all the orbitals of $h_j$
are contained in orbitals of $h_{j+1}$ for any $j\in\N$). $k$ has the
following properties.

\be

\item The orbitals of $h_{N_0}$ which contain the terminal orbitals of
$h_0$ are not contained in the orbitals of $k$.

\item No orbital of $h_{N_0}$ which is also an orbital of $h_{N_1+0}$
is also an orbital of $k$ (these are all terminal orbitals of
$h_{N_0}$ since $D$ is the result of using a growing tower operation).

\item All the non-terminal orbitals of $h_{N_0}$ are still properly
contained in the orbitals of $k$ since $k$ contains the non-terminal
orbitals of $h_{N_0+1}$.

\ee

Now replace $k$ and $h_{N_0}$ by sufficiently high powers of
themselves so that they satisfy the mutual efficiency condition and let
$h = [[h_{N_0},k], k]$.  The resulting $h$ has the following properties.

\be
\item $h$ has no orbitals intersecting the terminal orbitals of $h_0$.
\item $h$ has all the non-terminal orbitals of $h_{N_0}$
\ee

Now replace $h$ and $h_0$ by sufficiently high powers of themselves so
that the satisfy the mutual efficiency condition, and then replace
$h_0$ by $[[h_0, h],h]$.  Now replace $h_0$ and $h_{N_0}$ by
sufficiently high powers of themselves so they satisfy the mutual
efficiency condition.  Build the tower
\[
D' = \left\{(A_0,h_0)\right\}\cup \left\{(A_i,h_i)\,|\,i\geq N_0,
i\in\N\right\}.
\]

In this tower, $h_0$ admits only non-terminal orbitals, every
orbital of $h_0$ is properly contained in a non-terminal orbital of
$h_{N_0}$, and $h_0$ still has a copy of every non-terminal
orbital that it started with.  If we re-index the tower $D'$ and call
it $D$ again, then it satisfies all the old properties of the tower
$D$ found above, but its bottom element ($h_0$) has nice orbitals.  We
can now repeat this whole process for the sub-tower of $D$ starting
from level two and up, so that the new $h_1$ will admit all the
non-terminal orbitals that it started with, and other non-terminal
orbitals, and also will contain no terminal orbitals.  Inductively
proceed up the tower $D$, redefining all of the $h_i$, so that the new
tower $D$ satisfies the following properties.
\be

\item $A = \cup_{n\in\N} A_n$
\item For each $n\in\N$ with $n>0$, the orbitals of $h_{n-1}$ in $A$
are all contained inside the orbital $A_n$ of $h_n$.
\item For each $n\in\N$ with $n>0$, if $B$ is an orbital of $h_n$
which is not disjoint from the orbitals of $h_{n-1}$, then $B$
contains the closure of the union of the collection of orbitals of
$h_{n-1}$ that intersect $B$.

\ee 

\vspace{.1 in}

\underline{S6: Enabling wreath product structures by increasing fundamental domains}.

Now for each index $j\in\N$, inductively replace $h_j$ and $h_{j+1}$
by sufficiently high powers of themselves so that they satisfy the
mutual efficiency condition.   Note that each
signature (except $h_0$) may be replaced by progressively higher
powers of itself twice in this operation, but that once two signatures
are mutually efficient, replacing either signature by a higher power
of itself will still result in a pair that are mutually efficient.

Now any pair of adjacent signatures of the tower $D$ satisfy the
mutual efficiency condition.

\vspace{.1 in}

\underline{S7: Notes on dynamics with algebraic conclusions}.

For every $n\in\N$, define the subgroup $H_n=\langle
h_0,h_1,\ldots,h_n\rangle$ of $G$.  Now suppose that $n>0$.  Given any
two elements $f$, $g\in H_{n-1}$, since the supports of $f$ and $g$
are contained in the support of $h_{n-1}$, and since the support of
$h_{n-1}$ in any one orbital of $h_n$ is contained in a single
fundamental domain of $h_n$ in that orbital, we see that $f^{h_n^j}$
and $g^{h_n^k}$ have disjoint supports and therefore commute, whenever
$j\neq k$.  If $j=k$, then the product of the conjugated $f$ and $g$
is equal to the conjugate of the product of $f$ and $g$.  In
particular, the subgroup of $H_n$ consisting of finite products of
conjugates of elements of $H_{n-1}$ by $h_{n}$ is isomorphic to
$\sum_{j\in\Z}H_{n-1}$, where the indexing factor $j$ represents the
power of $h_{n}$ used in the conjugation of the element from $H_{n-1}$
under consideration.  But we can write any element of $H_{n}$ as a
product of an integer power of $h_{n}$ with a product of conjugates of
elements of $H_{n-1}$ by integer powers of $h_{n}$; in short,
$H_{n}\cong H_{n-1}\wr\Z$, where the $\Z$ factor is the subgroup
$\langle h_n\rangle$ of $H_n$.

Now, $H_0 \cong \Z$, so $H_1\cong \Z\wr\!\Z$, $H_2\cong
(\Z\!\wr\!\Z)\!\wr\!\Z$, and etc., so that $H_n\cong
((\cdots(\Z\!\wr\!\Z)\!\wr\!\Z)\cdots\wr\!\Z$ where the finite wreath product
has $n+1$ factors of $\Z$.  In particular, the ascending union
$H=\langle h_0,\,h_1,\,\ldots \rangle\cong(\Z\wr)^{\infty}$.

\qquad$\diamond$

\bl
\label{deepWreathZ}

If $G$ is a subgroup of $\ploi$ and $G$ admits a deep tower, then $G$
has a subgroup of the form $^{\infty}(\wr\Z)$.

\el

pf:

We will use a similar technique to the proof of Lemma
\ref{tallZWreath}, although the analysis in this case is much simpler.

If $G$ admits a transition chain of length two, then by Theorem
\ref{tChainB}, $G$ admits an embedded copy of $B$, and $B$ contains
copies of $^{\infty}(\wr\Z)$, so let us assume that $G$ admits no
transition chains of length two.

Since $G$ is admits no transition chains of length two, any tower for
$G$ is exemplary.  In particular, let $E =
\left\{(A_{-i},g_{-i})\,|\,i\in\N\backslash{0}\right\}$ be an
exemplary deep tower for $G$ where the indexing respects the order on
the elements of the tower.  Improve $E$ by replacing the signatures of
$E$ with sufficiently high powers of themselves so that given any
negative integer $i$, then $g_{i-1}$ and $g_i$ satisfy the mutual
efficiency condition.

Let $A=A_{-1}=(a,b)$.  Since $E$ is exemplary, we see that $A$ is
actually an orbital of the subgroup $H\leq G$, where $H =
\left<S_E\right>$.  For all $i\in\N$ with $i>1$, inductively improve
$E$ (induct on increasing $i\in\N$ in the following discussion) by
replacing the signatures of $E$ according to the following three step
process.  

First, let $h_{-i} = [[g_{-i},g_{-i+1}], g_{-i+1}]$.  

Second, define the new $g_{-i}$ to be $h_{-i}$.  

Third, replace the elements $g_{-i+1}$, $g_{-i}$, and $g_{-i-1}$
with sufficiently high powers of themselves, so that given any index
$j\in\N$, the elements $g_{-j}$ and $g_{-j-1}$ satisfy the mutual
efficiency condition (observe that if $i>3$, then $g_{-i+1}$ and
$g_{-i+2}$ will now still satisfy the mutually efficiency condition,
since we are only replacing $g_{-i+1}$ by higher powers of itself, and
these two signatures were already mutually efficient, a similar
argument shows that $g_{-i-1}$ and $g_{-i-2}$ will be mutually
efficient after this operation as well).  

Since $A_{-i}\subsetneq A_{-i+1}$ for all integers $i>1$, we see that the
resultant set of signed orbitals is still a tower (and with the same
order), so that this inductive definition will simply improve our
tower $E$.  Observe further that given any $k\in\N$, then the orbitals
of $g_{-k-1}$ are all properly contained in the orbitals of $g_{-k}$.

Define the set $\Gamma_i=\left\{g_j\,|\,j\leq i, j\in\Z\right\}$ for
each negative integer $i$.  For each negative integer $i$, define
$H_i=\langle \Gamma_i\rangle$.  For such $i$, the orbitals of $H_i$
are actually the orbitals of $h_i$, since all orbitals of the elements
$g_{k}$ with $k<i$ are contained in the orbitals of $g_i$.
Furthermore, for any such $i<-1$, the orbitals of $g_i$ are contained
in the orbitals of $g_{i+1}$ in such a way that in any individual
orbital $B$ of $g_{i+1}$, the support of $g_{i}$ in $B$ is actually
fully contained inside a single fundamental domain of $g_{i+1}$ on
$B$.  In particular, $H_i\cong H_{i-1}\!\wr\!\Z$, where the $\Z$
factor comes from the subgroup $\langle g_i\rangle$ of $H_i$.  But now
inductively, since each generator generates a group isomorphic to
$\Z$, we see that $H_1\cong\,^{\infty}(\wr\Z)$.

\qquad$\diamond$

\bl
\label{biWreathZ}

If $G$ is a subgroup of $\ploi$ and $G$ admits a bi-infinite tower, then $G$
has a subgroup of the form $(\wr\Z\wr)^{\infty}$.

\el

pf:

This follows immediately from the previous two lemmas, where first one
improves the non-negative tower, and then one improves the negative
tower (using the element with index $0$ as the top element).

\qquad$\diamond$

\subsection{$W$ in arbitrary non-solvable subgroups of $\ploi$}
We are now ready to complete the proof of Theorem
\ref{NSC}, restated below.

\bt

Let $H$ be a subgroup of $\ploi$.  $H$ is non-solvable if and only if
$W$ embeds in $H$.
\et

pf:

Since $W$ is non-solvable, any group which contains an embedded copy
of $W$ will be non-solvable as well, therefore we need only show that
if $H$ is non-solvable then $H$ contains a copy of $W$.

Suppose therefore that $H$ is a non-solvable subgroup of $\ploi$.  By
Lemma \ref{WInStuff} we know that $W$ embeds in $^{\infty}(\wr\Z)$ and
$(\Z\wr)^{\infty}$, and therefore also into $(\wr\Z\wr)^{\infty}$.
Therefore, if $H$ admits infinite towers then we already have the
result, so let us assume that $H$ does not admit infinite towers.  In
this case, we further have that $H$ admits towers of arbitrary finite
height, $H$ is balanced, and $H$ admits no transition chains of length
two.

Since $H$ does not admit infinite towers, the depth of any signed
orbital of $H$ is well defined and finite.  Since $H$ is not the
trivial group, $H$ has a non-empty collection of orbitals.  The
analysis now breaks into two cases.

\underline{{\bf Case 1:}}

Suppose $H$ admits no orbital that supports towers of arbitrary
height.  In this case the depth of any orbital of $H$ is well defined,
every orbital of $H$ has finite depth, and given any $n\in\N$, $H$ has
orbitals with depth greater than $n$.

Now, pick an element $g^1_1$ of $H$ so that $\hat{T}_1 =
\left\{(B_1^1,g_1^1)\right\}$ is a tower of height one for $H$.
$g^1_1$ will be our generator for $W_1$.  $g_1^1$ has finitely many
orbitals, and so there is a maximum depth $j_1$ of the orbitals of $H$
that are not disjoint from the support of $g_1^1$.  We will now pick
our remaining generators from the group $H^{(j_1)}$, the $j_1$'st
derived subgroup of $H$.  We note that no element in $H^{(j_1)}$ can
have support intersecting $g_1^1$, since $H^{(j_1)}$ has trivial
support over the orbitals of $H$ of depth less than or equal to $j_1$
as a consequence of the details of the proof of point 4 of Lemma
\ref{productOrbitals}. (We state the main idea of that proof, which is
in \cite{bpgsc}.  If $\iota$ is an element of $H$ and $(Z,\iota)$ is a
depth one signed orbital of $H$, then we note that there is no element
of $H$ with an orbital containing $\overline{Z}$, the closure of $Z$.
In particular, $Z$ cannot be the orbital of any commutator of elements
of $H$, nor of a finite product of commutators.  Hence, each time we
pass to a commutator subgroup, we lose all elements which support the
depth one orbitals of the original group.)

We also observe that $H^{(j_1)}$ still admits towers of arbitrary
height, and infinitely many orbitals, of arbitrary finite depth.  We
now find a tower $\hat{T}_2 =
\left\{(B_1^2,g_1^2),(B_2^2,g_2^2)\right\}$ for $H^{(j_1)}\leq H$ of
height two.  Now the signatures of $\hat{T}_2$ admit a finite total
number of orbitals, and therefore the union of this collection of
element orbitals is contained in the union of the collection of
orbitals of $H$ of depth less than some integer $j_2>j_1$.  We
therefore will pick a tower $\hat{T}_3=\left\{(B_1^3,g_1^3),
(B_2^3,g_2^3), (B_3^3,g_3^3)\right\}$ for $H^{(j_2)}$ which has
signatures whose supports must be disjoint from the supports of the
signatures of the first two towers $\hat{T}_1$ and $\hat{T}_2$ (by
using elements from $H^{(j_2)}$, for example).  We can continue in
this fashion to inductively define towers $\hat{T}_k$ and integers
$j_{k-1}$ for each positive integer $k$ so that the integers $j_{k}$
are always getting larger, and so that the towers $\hat{T}_k$ always
have height $k$ and have signatures which are disjoint in support from
the signatures of the previous towers.

Let $k\in\N$.  Let $\hat{G}_k$ represent the group generated by the
signatures of $\hat{T}_k$.  We can use the techniques of the proof of
Lemma \ref{deepWreathZ} to replace $\hat{T}_k$ with a new tower $T_k$
supported by the subset of the orbitals of $H$ that support $\hat{T}_k$, so
that the signatures of $T_k$ generate a group $G_k$ isomorphic to
$W_k$.  Do this for all $k\in\N$.

Now the union of all the signatures of all of the towers $T_k$ forms a
collection of generators of a group isomorphic to $W$.

\underline{{\bf Case 2:}}

Suppose now that $H$ admits an orbital $A$ that supports towers of
arbitrary height.

If $A$ is not an orbital of any element of $H$ then $A$ can be written
as a union of an infinite collection of nested element orbitals of
$H$, so that $H$ would then admit an infinite tower, therefore there
is an element $d$ of $H$ so that $(A,d)$ is a signed orbital of depth
one for $H$.

We will now restrict our attention to a special subgroup $H_d$ of $H$
which is directed by the element $d$, in a sense that will be made
clear.  Given any element $h\in H$, let $k_h$ and $j_h$ represent the
smallest positive integers so that $h^{k_h}$ and $d^{j_h}$ satisfy the
mutual efficiency condition.  Let
\[
\Gamma_d = \left\{[[h^{k_h},d^{j_h}], d^{j_h}]\,|\, h\in H\right\}\cup\left\{d\right\}.
\]
The elements of $\Gamma_d$ have all of their orbitals properly
contained inside the orbitals of $d$, and since the orbital $A$ of $H$
admits towers of arbitrary height, and any element orbital $B$ which
is properly contained inside $A$ will be realized as an orbital of
some element $g$ of $\Gamma_d$ (note that it does not matter that we
passed to high powers to guarantee the mutual efficiency condition),
we see that the group $H_d = \langle \Gamma_d\rangle$ admits towers of
arbitrary height. We now observe that given any finite set $X$ of
elements of $H_d$ that individually do not support any signed orbitals
of depth one for $H_d$, and a finite tower $T$ for $H_d$ which also
contains no signed orbital of depth one, we can find a minimal power
$k_X$ of $d$ so that so that the tower $T^{d^{k_X}}$ for $H_d$ induced
from $T$ via conjugation of the signatures of $T$ by $d^{k_X}$ will
have all of its signatures with disjoint support from the signatures
of $X$, since we can conjugate the tower to be arbitrarily near to an
end of an orbital of $d$.  Therefore, for each positive integer $n$,
let $\tilde{T}_n$ be a tower for $H_d$ of height $n$.  Now inductively
define towers $\hat{T}_n$ which are towers induced from the
$\tilde{T}_n$ by conjugation by powers of $d$ so that given any
positive integer $k$, the tower $\hat{T}_k$ has signatures whose
supports are all disjoint from the signatures of the towers
$\hat{T}_j$ whenever $j<k$ is a positive integer.

Now again apply the techniques of the proof of Lemma
\ref{deepWreathZ} to improve the towers $\hat{T}_n$ to new towers $T_n$
so that for each positive integer $n$, the signatures of the tower
$T_n$ generate a group isomorphic with $W_n$, while preserving the
conditions that the signatures of distinct towers $T_k$ and $T_j$ have
disjoint supports from each other, so that the union of all of the
signatures of all of the towers $T_k$ is a set of generators of a
group isomorphic with $W$ in $H$.

\qquad$\diamond$

\newpage
\bibliographystyle{amsplain}

\bibliography{ploiBib}

\providecommand{\bysame}{\leavevmode\hbox to3em{\hrulefill}\thinspace}
\providecommand{\MR}{\relax\ifhmode\unskip\space\fi MR }
\providecommand{\MRhref}[2]{%
  \href{http://www.ams.org/mathscinet-getitem?mr=#1}{#2}
}
\providecommand{\href}[2]{#2}
\begin{thebibliography}{10}

\bibitem{bpasc}
Collin Bleak, \emph{An algebraic classification of some solvable groups of
  homeomorphisms}, preprint (2005), 1--29.

\bibitem{bpgsc}
\bysame, \emph{A geometric classification of some solvable groups of
  homeomorphisms}, preprint (2005), 1--19.

\bibitem{BDiss}
\bysame, \emph{Solvability in groups of piecewise-linear homeomorphisms of the
  unit interval}, Binghamton University, 2005, Dissertation.

\bibitem{BrinU}
Matthew~G. Brin, \emph{The ubiquity of {T}hompson's group {$F$} in groups of
  piecewise linear homeomorphisms of the unit interval}, J. London Math. Soc.
  (2) \textbf{60} (1999), no.~2, 449--460.

\bibitem{BrinEG}
\bysame, \emph{Elementary amenable subgroups of {R}. {T}hompson's group {$F$}},
  Internat. J. Algebra Comput. \textbf{15} (2005), no.~4, 619--642.
  \MR{MR2160570}

\bibitem{BGAutomorphisms}
Matthew~G. Brin and Fernando Guzm{\'a}n, \emph{Automorphisms of generalized
  {T}hompson groups}, J. Algebra \textbf{203} (1998), no.~1, 285--348.

\bibitem{BSPLR}
Matthew~G. Brin and Craig~C. Squier, \emph{Groups of piecewise linear
  homeomorphisms of the real line}, Invent. Math. \textbf{79} (1985), no.~3,
  485--498.

\bibitem{picric}
\bysame, \emph{Presentations, conjugacy, roots, and centralizers in groups of
  piecewise linear homeomorphisms of the real line}, Comm. Algebra \textbf{29}
  (2001), no.~10, 4557--4596.

\bibitem{BrownFinite}
Kenneth~S. Brown, \emph{Finiteness properties of groups}, Journal of Pure and
  Applied Algebra \textbf{44} (1987), 45--75.

\bibitem{BCS}
J.~Burillo, S.~Cleary, and M.~I. Stein, \emph{Metrics and embeddings of
  generalizations of {T}hompson's group {$F$}}, Trans. Amer. Math. Soc.
  \textbf{353} (2001), no.~4, 1677--1689 (electronic).

\bibitem{SteinPLGroups}
Melanie Stein, \emph{Groups of piecewise linear homeomorphisms}, Trans. Amer.
  Math. Soc. \textbf{332} (1992), no.~2, 477--514. \MR{MR1094555 (92k:20075)}

\end{thebibliography}
\end{document}